# THE TRIANGULAR THEOREM OF THE PRIMES: BINARY QUADRATIC FORMS AND PRIMITIVE PYTHAGOREAN TRIPLES

J. A. PEREZ

ABSTRACT. This article reports the occurrence of binary quadratic forms in primitive Pythagorean triangles and their geometric interpretation. In addition to the well-known fact that the hypotenuse, $z$, of a right triangle, with sides of integral (relatively prime) length, can be expressed as the sum of two squares, $z = a^2 + b^2$ (where $a, b \in \mathbb{N}$, such that $a > b > 0$, $a \not\equiv b \pmod 2$ and $\gcd(a,b) = 1$), it is shown that the sum of the two sides, $x$ and $y$, can also be expressed as a binary quadratic form, $x + y = (a+b)^2 - 2b^2$. Similarly, when the radius of the inscribed circle is taken into account, $r = b \cdot (a-b)$, a third binary quadratic form is found, namely $(x+y) - 4r = z - 2r = (a-b)^2 + 2b^2$. The three quadratic representations accommodate positive integers whose factorizations can only include primes $p$ represented by the same type of binary quadratic forms, i.e. $p \equiv 1, 5 \pmod 8$, $p \equiv 1, 7 \pmod 8$, and $p \equiv 1, 3 \pmod 8$, respectively. For all three types of binary quadratic forms, when the positive integers represented are prime, such representations are unique. This implies that all odd primes can be geometrically incorporated into primitive Pythagorean triangles.

## 1. Introduction

The representation of natural numbers as sums of integer squares has been a problem studied by number theorists for centuries [9,10,11,22,29]. Of particular interest are the cases in which the natural numbers represented in this way are prime. Fermat was the first to indicate that for an odd prime $p$,

$$p = a^2 + b^2, \ a, b \in \mathbb{N} \ \Leftrightarrow \ p \equiv 1 \pmod 4, \tag{1.1}$$

as eventually proved by Euler [9,11]. Similarly, Fermat showed that

$$p = a^2 + 2b^2, \ a, b \in \mathbb{N} \ \Leftrightarrow \ p \equiv 1, 3 \pmod 8, \tag{1.2}$$

$$p = a^2 + 3b^2, \ a, b \in \mathbb{N} \ \Leftrightarrow \ p = 3 \ \text{or} \ p \equiv 1 \pmod 3, \tag{1.3}$$

and again Euler provided proofs for both statements [9,11]. The consequent development of quadratic reciprocity and the theory of quadratic forms extended these initial results considerably [9,10,13,22]. Further progress was achieved with the advent of higher reciprocity, class field theory and complex multiplication [9,11,18,27], as well as many other aspects of algebraic and analytical number theory [4,21,30].

Equation (1.1) merits particular attention because of its direct association with the Pythagorean equation $x^2+y^2=z^2$. It can be said that, for primes $p\equiv 1 \pmod 4$, there are integers $a$ and $b$ such that

$$(2ab)^2+(a^2-b^2)^2=(a^2+b^2)^2=p^2\ , \tag{1.4}$$

therefore making $p$ the hypotenuse of a right triangle with sides $a^2-b^2$, $2ab$, and $p$ [7,20,30]. Fermat called (1.4) "the fundamental theorem on right triangles" [30], clearly giving some importance to the geometric interpretation of (1.1).

These observations suggest that other groups of primes, represented by alternative binary quadratic forms, could also be geometrically expressed within Pythagorean triangles. This article describes such findings.

## 2. Primitive Pythagorean triples

Any set of integers $x,y,z\in\mathbb{N}$, that satisfy the Pythagorean equation $x^2+y^2=z^2$, is defined as a *Pythagorean triple*. Furthermore, the triple is said to be *primitive* if $\gcd(x,y,z)=1$ [7]. All the primitive triples can be generated by the well-known parametrization

$$x=2ab \qquad y=a^2-b^2 \qquad z=a^2+b^2 \tag{2.1}$$

for integers $a,b\in\mathbb{N}$, such that $a>b>0$, $\gcd(a,b)=1$, and $a\not\equiv b \pmod 2$ [7,20,29]. This can be used to define a *primitive Pythagorean triangle* as any right triangle whose sides are of integral length and relatively prime (see Figure 1). It is a geometric property of all Pythagorean triangles that the radius $r$ of their inscribed circle is always an integer [7,19]. Using the parametrization (2.1), it is simple to deduce that the value of $r$ for a primitive Pythagorean triangle is given by

$$r=b\cdot(a-b)\ . \tag{2.2}$$

This can be derived from a classical observation that relates the two sides $x,y$ and the hypotenuse $z$ of a right triangle to the radius $r$ of its inscribed circle [19],

$$z=x+y-2r\ , \tag{2.3}$$

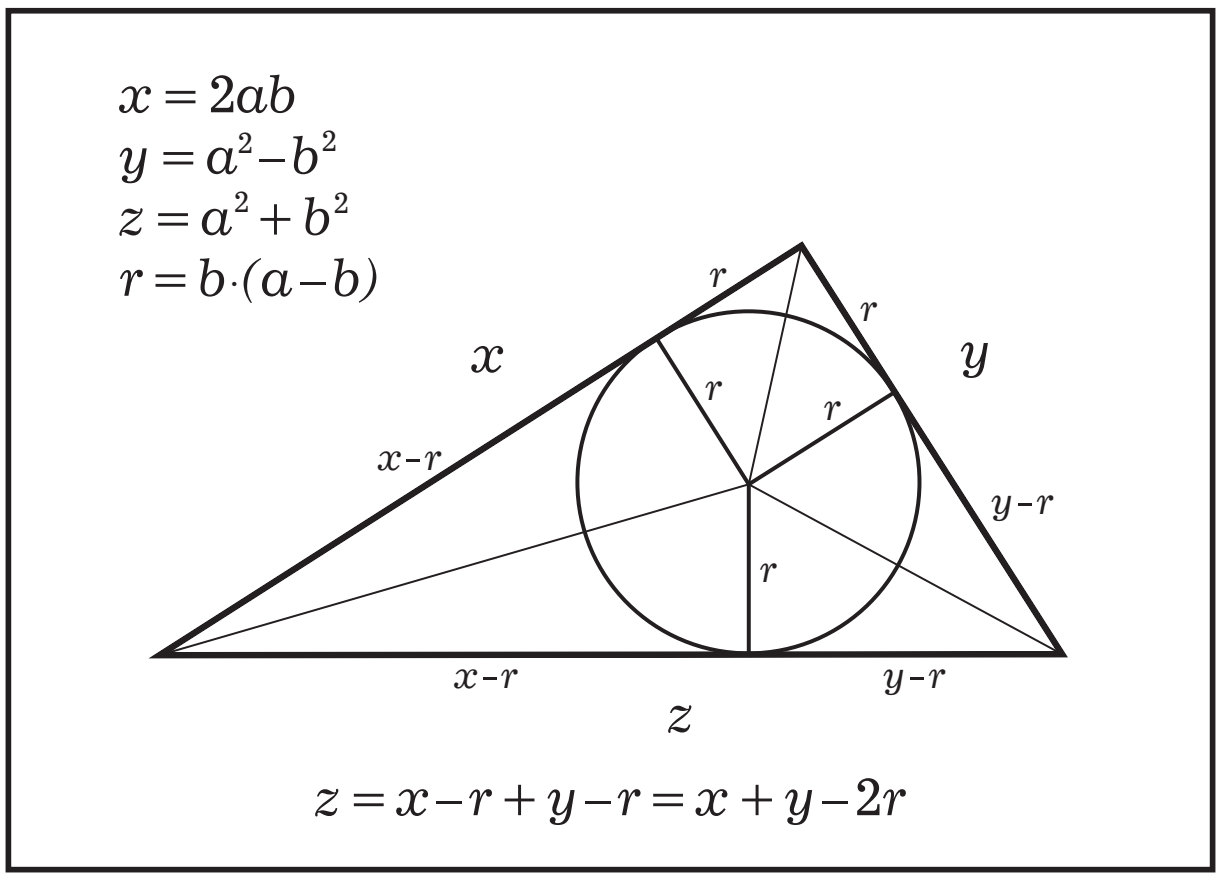

Figure 1. The parametrization of primitive Pythagorean triangles.

as easily deduced from Figure 1. The parametrization (2.1) and the relation (2.3) readily combine to deliver (2.2):

$$2r = x + y - z = 2ab + (a^2 - b^2) - (a^2 + b^2) = 2ab - 2b^2 = 2b\cdot(a-b)\ .$$

The geometric and/or algebraic observations described by (2.1)–(2.3) can be extended further. As illustrated in Figure 2, the parametrization (2.1) can be used to deduce that:

*a)* the sum of the two sides $x, y$ can be expressed as a binary quadratic form, i.e.

$$x + y = 2ab + (a^2 - b^2) = a^2 + 2ab - b^2 = (a+b)^2 - 2b^2\ ; \tag{2.4}$$

*b)* as established by (2.3) [19], the substraction of two radii from the sum of the two sides $x, y$ renders the hypothenuse $z$, expressed as a sum of two squares,

$$\begin{aligned} x - r &= 2ab - b\cdot(a-b) = ab + b^2\ ,\\ y - r &= a^2 - b^2 - b\cdot(a-b) = a^2 - ab\ ,\\ x + y - 2r &= z = a^2 + b^2\ ; \end{aligned} \tag{2.5}$$

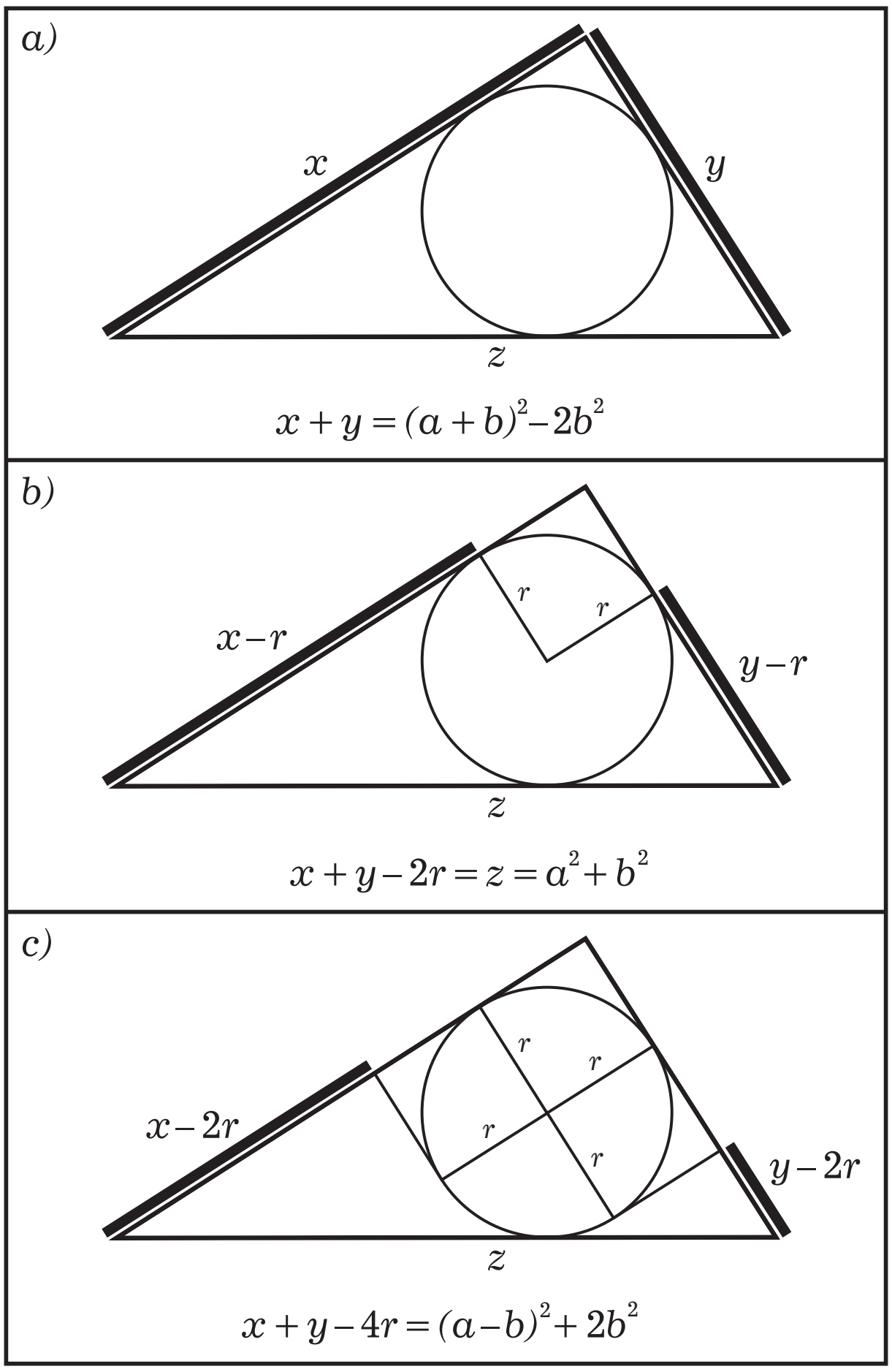

FIGURE 2. Binary quadratic forms in primitive Pythagorean triangles.

*c)* the substraction of four radii from the sum of the two sides $x, y$ generates another binary quadratic form,

$$\begin{aligned} x-2r &= 2ab-2b\cdot(a-b) = 2b^2\,, \\ y-2r &= a^2-b^2-2b\cdot(a-b) = a^2-2ab+b^2 = (a-b)^2\,, \\ x+y-4r &= (a-b)^2+2b^2\,. \end{aligned} \tag{2.6}$$

Hence, equations (2.4) – (2.6) represent three distinct binary quadratic forms that are all naturally associated with the geometry of primitive Pythagorean triangles (Figure 2).

## 3. Binary quadratic forms

The geometric role that the three binary quadratic forms (2.4)–(2.6) clearly display in all primitive Pythagorean triangles (Figure 2) prompts a detailed characterization of the integers represented by these forms, within the context of the parametrization (2.1).

3.1. $\boldsymbol{x+y-2r=z=a^2+b^2}$**.** As indicated in (1.4), the hypotenuse $z=a^2+b^2$ of a primitive Pythagorean triangle can be seen as representing, for some value of $a$ and $b$, an odd prime $p\equiv 1 \pmod 4$. A more complete statement is to affirm that any integer $N\in\mathbb{N}$, represented by the sum of two squares, $N=a^2+b^2$, such that $a,b\in\mathbb{N}$, $a>b>0$, $a\not\equiv b \pmod 2$ and $\gcd(a,b)=1$, can only be divided evenly by factors that are also written as sums of two squares. When these factors are prime, they are of the form $p\equiv 1 \pmod 4$ [7,10,11]. Requiring $a$ and $b$ to have opposite parities, $a\not\equiv b \pmod 2$, makes $N\equiv 1 \pmod 2$, thus excluding even divisors. Consequently, it can be established that the hypotenuses of the set of all primitive Pythagorean triangles, as previously defined, represent all possible integers $N$ whose factorizations include only primes of the form $p\equiv 1 \pmod 4$. Hence, define the set of integers

$$S_{1,5} := \Big\{ N_{1,5} = \prod_{j=1}^{n} p_j^{k_j} = a^2+b^2,\, a,b\in\mathbb{N},\, a>b>0, \\ \gcd(a,b)=1,\, a\not\equiv b \pmod 2,\, p_j\equiv 1,5 \pmod 8 \Big\}\,, \tag{3.1}$$

where $\prod_{j=1}^{n} p_j^{k_j}$ is the prime canonical factorization of $N_{1,5}$.

Since $p_j\equiv 1 \pmod 4 \equiv 1,5 \pmod 8$, it is always the case that $N_{1,5}\equiv 1,5 \pmod 8$, and this explains the nomenclature used in (3.1). That all possible integers $N_{1,5}$, as defined in (3.1), can be written as sums of two squares, and that at least one of such representations is primitive, are the consequence of three established facts. First, the well-known identity

$$(a^2+b^2)(c^2+d^2) = (ac\mp bd)^2 + (ad\pm bc)^2 \tag{3.2}$$

implies that, if two numbers are sums of two squares, then their product is also a sum of two squares [9,11,22]. Secondly, as originally proved by Euler using an infinite descent, when $a$ and $b$ are relatively prime, every factor of $N_{1,5}=a^2+b^2$ is itself a sum of two squares [9,11]. Thirdly, for any integer $N_{1,5}=\prod_{j=1}^{n} p_j^{k_j}$, such that $p_j\equiv 1,5 \pmod 8$, the total number of proper (primitive) representations as a sum of two squares is given by $2^{n-1}$, where $n$ is the number of distinct odd prime divisors of $N_{1,5}$ (see, for example Section 8 of chapter VI in [10]; also [22]). As independently proved [7,22], when $N_{1,5} = p_j = a^2 + b^2$, the primitive representation is therefore unique. Accordingly, the number of primitive Pythagorean triangles whose hypotenuse represent an integer of $S_{1,5}$ is dependent on the factorization of $N_{1,5}$.

3.2. $\boldsymbol{x + y = a^2 + 2ab - b^2 = (a + b)^2 - 2b^2}$**.** Equation (2.4) indicates that the sum of the two sides of a primitive Pythagorean triangle can be represented by a binary quadratic form, in a way that parallels the case of the hypotenuse $z = a^2 + b^2$. This new form $F(a,b) = (a+b)^2 - 2b^2$ can be found to be equivalent to another form $G(a',b')$ according to

$$x + y = a^2 + 2ab - b^2 = (a+b)^2 - 2b^2 = a'^2 - 2b'^2 \,, \tag{3.3}$$

that results from the transformation

$M = \begin{pmatrix} 1 & 0 \\ -1 & 1 \end{pmatrix}$, which induces the change $(a+b,\, b) \to (a,b)$.

The transformation matrix $M$ has determinant $\det(M) = +1$; this implies that $M$ is unimodular and justifies the equivalence $F \sim G$ [27,30].

The equivalence $F \sim G$ proves useful when trying to determine what integers can be represented by $F(a,b) = a^2 + 2ab - b^2 = (a+b)^2 - 2b^2$. As is well-documented in the literature (see [4,10,11]; also exercise 2.23 in [9]), all primes $p \equiv 1,7 \pmod 8$ can be represented by the quadratic form $G(a',b') = a'^2 - 2b'^2$. Therefore, by making use of (2.4), it is possible to establish that these primes could, in principle, be represented by $p = (a+b)^2 - 2b^2 = x + y$, where $x$ and $y$ are the sides of a primitive Pythagorean triangle. It is trivial to see that, for $G(a',b') = a'^2 - 2b'^2$ to represent an odd prime, as well as requiring $\gcd(a',b') = 1$, it is necessary that $a' \equiv 1 \pmod 2$. Since it is always true for any primitive Pythagorean triangle that $\gcd(a,b) = 1$, and $a \not\equiv b \pmod 2$, it is elementary to establish that $\gcd(a+b,b) = 1$ and $a + b \equiv 1 \pmod 2$, thus confirming that the form $F(a,b) = (a+b)^2 - 2b^2$ is capable of representing odd primes.

It should be pointed out that the discriminant of the forms $F$ and $G$ is positive ($\Delta = 8$), and that makes them indefinite. This would imply that, in principle, they could take both positive and negative values [10,22,27]. However, the possibility of negative values is eliminated by the condition $a > b > 0$, as required by the parametrization (2.1) [7,20,29]. Consequently, the general principles of the theory of quadratic forms, as applicable to positive definite binary forms [10,27], should also be applicable to the forms $F$ and $G$, as described in (3.3), provided that $a > b > 0$. The crucial role played by this condition is illustrated by the fact that the quadratic form $G(a',b') = a'^2 - 2b'^2$, which does not impose any restrictions on the relative values of $a'$ and $b'$, generates infinite representations for primes $p \equiv 1,7 \pmod 8$ [22,27]. As example, take the four smallest primes $p \equiv 1,7 \pmod 8$, and examine their first five representations, i.e. those using the smallest possible (positive) values for $a$ and $b$:

$$7 = (2+1)^2 - 2\cdot 1^2 = (2+3)^2 - 2\cdot 3^2 = (4+9)^2 - 2\cdot 9^2 = (8+19)^2 - 2\cdot 19^2 = (22+53)^2 - 2\cdot 53^2$$

$17 = (3+2)^2 - 2\cdot2^2 = (3+4)^2 - 2\cdot4^2 = (7+16)^2 - 2\cdot16^2 = (11+26)^2 - 2\cdot26^2 = (39+94)^2 - 2\cdot94^2$

$23 = (4+1)^2 - 2\cdot1^2 = (4+7)^2 - 2\cdot7^2 = (6+13)^2 - 2\cdot13^2 = (18+43)^2 - 2\cdot43^2 = (32+77)^2 - 2\cdot77^2$

$31 = (4+3)^2 - 2\cdot3^2 = (4+5)^2 - 2\cdot5^2 = (10+23)^2 - 2\cdot23^2 = (14+33)^2 - 2\cdot33^2 = (56+135)^2 - 2\cdot135^2$

In all four cases, the first representation is the only one that meets the condition $a > b$. This observation might be taken to suggest two things: *i)* first, all primes $p \equiv 1,7 \pmod 8$, expressed by the form $F(a,b) = (a+b)^2 - 2b^2$, have at least one representation where $a > b > 0$; *ii)* secondly, such a representation is unique. The following two theorems provide the proofs.

**Theorem 3.1.** *Any prime number* $p \equiv 1,7 \pmod 8$*, expressed by the quadratic form* $F(a,b) = (a+b)^2 - 2b^2$*, where* $a, b \in \mathbb{N}$, $a \not\equiv b \pmod 2$ and $\gcd(a,b) = 1$, *has at least one representation where* $a > b > 0$.

*Proof.* By contradiction. In order to reach this conclusion, assume the opposite is true, i.e. that for *all* infinite representations of $p \equiv 1,7 \pmod 8$ using the binary quadratic form $F(a,b) = (a+b)^2 - 2b^2$, $b > a$ whenever $b > 0$.

As a preliminary fact, it is important to verify that, when $b > 0$, it is also true that $a > 0$. Since $p = F(a,b) > 0$, it is the case that $(a+b)^2 > 2b^2$, which in turn means that $|a+b| > +\sqrt{2}b$. Consequently, for a fixed value of $b$, there will always be two possible values of $a$. Considering the first possibility, $a_1 + b > 0$, the implication is that $a_1 + b > \sqrt{2}b$, and $a_1 > (\sqrt{2}-1)\,b > 0$, as required. The second possibility, $a_2 + b < 0$, implies that $a_2 < 0$, since $a_2 + b = -(a_1 + b)$ and

$$a_2 = -(a_1 + 2b) < 0. \tag{3.4}$$

It is also important to observe that, in all cases, $a \neq b$. This can be deduced by contradiction since, if $a = b$, then

$$F(a,b) = (a+b)^2 - 2b^2 = (2a)^2 - 2a^2 = 2a^2, \tag{3.5}$$

which is incompatible with $p$ being an odd prime, $p \equiv 1 \pmod 2$.

Applying the Well-Ordering Principle [7] to the infinite set of representations of $p = (a+b)^2 - 2b^2$, where $b > 0$, consider the representation with the smallest possible value of $b$: $p = (a_1 + b_1)^2 - 2b_1^2$, where $b_1 > a_1$, according to the starting assumption. Consider also a second representation of the set $p = (c_1 + d_1)^2 - 2d_1^2$, such that $c_1 = a_1$. Therefore,

$$p = (a_1 + b_1)^2 - 2b_1^2 = (c_1 + d_1)^2 - 2d_1^2 = (a_1 + d_1)^2 - 2d_1^2, \tag{3.6}$$

and

$$a_1^2 + 2a_1b_1 - b_1^2 = a_1^2 + 2a_1d_1 - d_1^2,$$

$$2a_1(b_1 - d_1) = b_1^2 - d_1^2 = (b_1 + d_1)(b_1 - d_1),$$

$$d_1 = 2a_1 - b_1. \tag{3.7}$$

The simplification leading to equation (3.7) can only be legitimately obtained if $d_1 \neq b_1$, but this is also compatible with $a_1 \neq b_1$, as required to avoid (3.5). The implication is that condition $c_1 = a_1$ can always be satisfied, i.e. $p$ always has two representations with identical values $c_1 = a_1$ but different values $d_1 \neq b_1$. As given by (3.7), all three possibilities regarding value and sign of $d_1$ must be considered:

(1) $d_1 = 2a_1 - b_1 = 0$
This possibility is incompatible with $p$ being prime, since it implies that $p = (c_1 + d_1)^2 - 2d_1^2 = c_1^2$ is a perfect square.

(2) $d_1 = 2a_1 - b_1 > 0$
As initially assumed, $b_1 > a_1$. Therefore, $d_1 = 2a_1 - b_1 = a_1 + (a_1 - b_1) < a_1$ and $d_1 < c_1$, hence contradicting the original assumption.

(3) $d_1 = 2a_1 - b_1 < 0$
It can then be stated that $b_1 - 2a_1 > 0$. Accordingly,

$$\begin{aligned} p &= (c_1 + d_1)^2 - 2d_1^2 = [a_1 + (2a_1 - b_1)]^2 - 2(2a_1 - b_1)^2 = [a_1 - (b_1 - 2a_1)]^2 - 2(b_1 - 2a_1)^2 \\ p &= [a_1 - 2(b_1 - 2a_1) + (b_1 - 2a_1)]^2 - 2(b_1 - 2a_1)^2 = (c_2 + d_2)^2 - 2d_2^2 , \end{aligned} \tag{3.8}$$

where the new representation has parameters

$$\begin{aligned} c_2 &= a_1 - 2(b_1 - 2a_1) = 5a_1 - 2b_1 , \\ d_2 &= b_1 - 2a_1 > 0 . \end{aligned} \tag{3.9}$$

It is important to determine the sign of $c_2$. The application of (3.4) to (3.9) implies that, for a fixed $d_2$, the alternative to $c_2$ takes the value

$$-(c_2 + 2d_2) = -(5a_1 - 2b_1 + 2b_1 - 4a_1) = -a_1 < 0 . \tag{3.10}$$

Therefore, it can be concluded that $c_2 = 5a_1 - 2b_1 > 0$. Since $b_1 > 2a_1$, it can also be deduced that

$$\begin{aligned} c_2 &= 5a_1 - 2b_1 < 5a_1 - 2(2a_1) = a_1 , \\ d_2 &= b_1 - 2a_1 < b_1 . \end{aligned} \tag{3.11}$$

As shown before, $c_2 \neq d_2$; hence, there are two possibilities:

(a) $c_2 > d_2$, contradicting the original assumption.

(b) $c_2 < d_2$, contradicting that the representation with $a_1$ and $b_1$ has the smallest possible value of $b$, while $b_1 < a_1$.

The overall conclusion is finally reached that for each prime $p \equiv 1,7 \pmod 8$ there is at least one representation $p = (a + b)^2 - 2b^2$ where $a > b > 0$. □

In the proof of the next theorem, use is made of the well-known formula [9,11]

$$(a^2 + kb^2)(c^2 + kd^2) = (ac \mp kbd)^2 + k(ad \pm bc)^2 , \qquad k = \pm 1, 2, 3, ... \tag{3.12}$$

that, when making $k = -2$ and adopting forms like $F(a,b)$, leads to the following lemma.

**Lemma 3.2.** *Consider two integers, p and q (prime or not), expressed by the binary quadratic form* $F(a,b) = (a + b)^2 - 2b^2$*:*

$$p = (a + b)^2 - 2b^2 \qquad\qquad q = (c + d)^2 - 2d^2 \tag{3.13}$$

*Their product* $N = p \cdot q$ *will be equally expressed by the form* $F(a,b) = (a + b)^2 - 2b^2$*, according to*

$$N = p\cdot q = (A+B)^2 - 2B^2 = (C+D)^2 - 2D^2 \tag{3.14}$$

*where*

$$A = ac + bd, \quad B = b(c+d) + d(a+b), \quad C = a(c+d) + d(a-b), \quad D = bc - ad. \tag{3.15}$$

*Proof.* Combining (3.12) and (3.13), it can be written that

$$\begin{aligned} p\cdot q &= \left[(a+b)^2 - 2b^2\right]\left[(c+d)^2 - 2d^2\right] = \left[(a+b)(c+d) + 2bd\right]^2 - 2\left[b(c+d) + (a+b)d\right]^2 \\ p\cdot q &= \left[(a+b)^2 - 2b^2\right]\left[(c+d)^2 - 2d^2\right] = \left[(a+b)(c+d) - 2bd\right]^2 - 2\left[b(c+d) - (a+b)d\right]^2. \end{aligned} \tag{3.16}$$

By making use of (3.14), the different terms in (3.16) can be expanded as follows:

$$\begin{aligned} &B = b(c+d) + d(a+b) = bc + ad + 2bd \\ &A + B = (a+b)(c+d) + 2bd = ac + bd + (bc + ad + 2bd) = ac + bd + B \\ &A = ac + bd \end{aligned}$$

and

$$\begin{aligned} &D = b(c+d) - d(a+b) = bc + bd - ad - bd = bc - ad \\ &C + D = (a+b)(c+d) - 2bd = a(c+d) + bc - bd = a(c+d) + bc - bd \pm ad \\ &C + D = a(c+d) + ad - bd + (bc - ad) = a(c+d) + d(a-b) + D \\ &C = a(c+d) + d(a-b) \end{aligned}$$

thus confirming the expressions given by (3.15). □

**Theorem 3.3.** *The representation of a prime number* $p \equiv 1,7 \pmod 8$ *by the binary quadratic form* $F(a,b) = (a+b)^2 - 2b^2$*, where* $a, b \in \mathbb{N}$*,* $a \not\equiv b \pmod 2$ and $\gcd(a,b) = 1$*, is always unique, provided that* $a > b > 0$.

*Proof.* By contradiction. To establish that the representation is unique, suppose the opposite is true, that is

$$p = (a+b)^2 - 2b^2 = (c+d)^2 - 2d^2, \tag{3.17}$$

where $a, b, c, d$ are all positive integers that satisfy the conditions $a > b$ and $c > d$. Using (3.17), it can be established that

$$\begin{aligned} &p\left[(a+b)^2 - (c+d)^2\right] = (a+b)^2\left[(c+d)^2 - 2d^2\right] - (c+d)^2\left[(a+b)^2 - 2b^2\right] \\ &p\left[(a+b)^2 - (c+d)^2\right] = 2\left[b^2(c+d)^2 - (a+b)^2 d^2\right]. \end{aligned} \tag{3.18}$$

Since $p$ is an odd prime, (3.18) implies that

$$\left[b(c+d) + (a+b)d\right]\left[b(c+d) - (a+b)d\right] \equiv 0 \pmod p. \tag{3.19}$$

Therefore, the existence of at least two representations implies compliance with congruence (3.19), which in turn requires that, either $b(c+d) + (a+b)d \equiv 0 \pmod p$, or $b(c+d) - (a+b)d \equiv 0 \pmod p$. The two possibilities need to be analysed separately.

(1) $b(c+d) + (a+b)d \equiv 0 \pmod p$
Since $b(c+d) + (a+b)d > 0$, the above congruence implies that there is a positive integer $k$ such that

$$b(c+d) + (a+b)d = kp. \tag{3.20}$$

(3.16), (3.17) and (3.20) can be combined to deduce that

$$p^2=[(a+b)^2-2b^2][(c+d)^2-2d^2]=[(a+b)(c+d)+2bd]^2-2[b(c+d)+(a+b)d]^2$$
$$p^2=(ac+ad+bc+bd+2bd)^2-2k^2p^2=[ac+bd+b(c+d)+(a+b)d]^2-2k^2p^2$$
$$p^2=(ac+bd+kp)^2-2k^2p^2,$$

thus obtaining

(3.21) $$p^2(1+2k^2)=(ac+bd+kp)^2.$$

The analysis of equation (3.21) requires evaluation of $ac+bd$ in relation to $p$. From the initial conditions $a>b$ and $c>d$, it can be inferred that $b\le a-1$ and $d\le c-1$ for all possible cases. Therefore,

(3.22) $$ac+bd\le ac+(a-1)(c-1)=2ac-(a+c)+1.$$

It can be assumed, without detriment to the proof, that $a\ge c$. Applying this condition to (3.22) generates

(3.23) $$ac+bd\le 2a^2-(a+c)+1<2a^2,$$

given that $a>1$, $c>1$ for all cases. Since $p=(a+b)^2-2b^2=a^2+2ab-b^2$, and $a>b$, it is elementary to establish that

(3.24) $$2p=2a^2+4ab-2b^2>2a^2.$$

(3.23) and (3.24) can then be used with (3.21) to produce

$$p^2(1+2k^2)=(ac+bd+kp)^2<(2p+kp)^2=p^2(2+k)^2$$

and establish the condition

(3.25) $$1+2k^2<(2+k)^2.$$

Condition (3.25) sets the boundaries for any possible values of $k$. It is elementary to find that the intersection points which equalize the two sides of (3.25) are given by the values $k=2\pm\sqrt{7}$ (see Figure 3). Consequently, (3.25) can only be satisfied within the range $2-\sqrt{7}<k<2+\sqrt{7}$. Since $k$ has to be a positive integer, the possible values to consider are $k=1, 2, 3$ or 4. These four possibilities can then be taken to (3.21) to evaluate the term $ac+bd$ as a function of $p$:

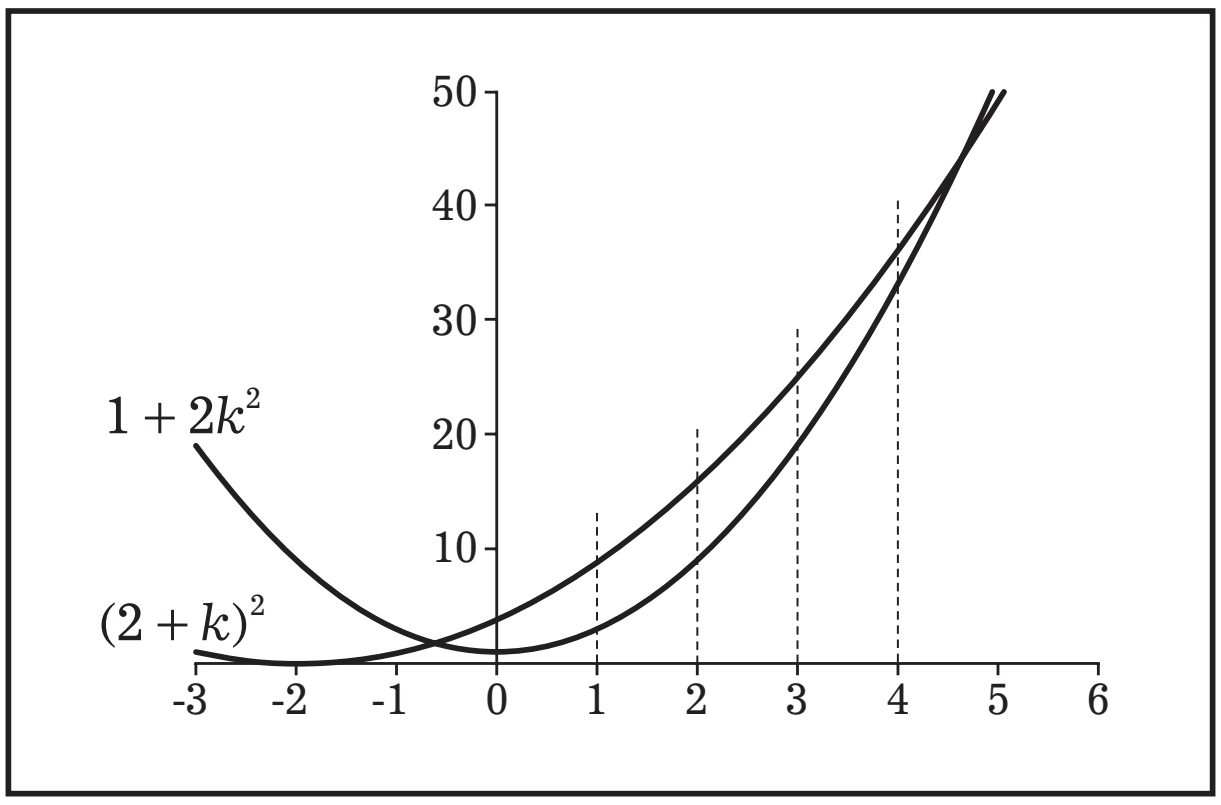

FIGURE 3. Analysis of (3.25): $1+2k^2<(2+k)^2$.

(a) $k=1$. Equation (3.21) becomes $3p^2=(t+p)^2$, where $t=ac+bd>0$. The positive solution is $t=\left(-1+\sqrt{3}\right)p>0$, which can never be an integer and rules out this possibility.

(b) $k=2$. (3.21) can be written as $9p^2=(t+2p)^2$, therefore $3p=t+2p$, with solution $t=p>0$, which satisfies condition (3.25).

(c) $k=3$. In this case, (3.21) becomes $19p^2=(t+3p)^2$. The corresponding positive solution is $t=\left(-3+\sqrt{19}\right)p>0$, which cannot be an integer and again rules out the possibility.

(d) $k=4$. For this last case, (3.21) can be written as $33p^2=(t+4p)^2$. The positive solution is $t=\left(-4+\sqrt{33}\right)p>0$, which once more cannot be an integer.

Since condition (3.25) can only be satisfied when $k=2$, making $ac+bd=p$, these two values can be taken to (3.20) to derive

$$2p=2ac+2bd=b(c+d)+(a+b)d=bc+ad+2bd$$

(3.26) $$2ac=bc+ad\,.$$

However, from the initial assumption, $b<a$ and $d<c$, implying that $bc<ac$ and $ad<ac$. Therefore, $bc+ad<ac+ac=2ac$, in direct contradiction of (3.26). It is to be concluded that the condition $b(c+d)+(a+b)d\equiv 0\,(\mathrm{mod}\,p)$ can never be met.

(2) $b(c+d)-(a+b)d\equiv 0\,(\mathrm{mod}\,p)$

The analysis of the above congruence requires assessment of the relative values of $b(c+d)$ and $(a+b)d$. Taking (3.17), it can be written that

(3.27) $$b=+\sqrt{\left[(a+b)^2-p\right]/2}$$ $$(c+d)=+\sqrt{p+2d^2}$$

and, therefore,

$$b(c+d)-(a+b)d=\sqrt{\left[(a+b)^2-p\right]/2}\,\sqrt{p+2d^2}-(a+b)d$$

$$b(c+d)-(a+b)d=\left(1/\sqrt{2}\right)\sqrt{(a+b)^2p+2(a+b)^2d^2-p^2-2d^2p}\;-(a+b)d$$

(3.28) $$b(c+d)-(a+b)d=\left(1/\sqrt{2}\right)\sqrt{2(a+b)^2d^2+p\left[(a+b)^2-2d^2\right]-p^2}\;-(a+b)d\,.$$

To evaluate equation (3.28), it will be sufficient to show that, in all cases, $3p>(a+b)^2-2d^2$. Again, from (3.17)

(3.29) $$3p=(a+b)^2-2d^2+(c+d)^2-2b^2+p\,,$$

where it is necessary to consider value and sign of the term $(c+d)^2-2b^2$. There are three possibilities:

(a) $(c+d)^2-2b^2=0$. This implies that $p=(c+d)^2-2d^2=2b^2-2d^2=2(b^2-d^2)$, incompatible with $p$ being an odd prime.

(b) $(c+d)^2-2b^2>0$. From (3.29), the implication is that $3p>(a+b)^2-2d^2$, as required.

(c) $(c+d)^2-2b^2<0$. Since $a>b$, then $p=(a+b)^2-2b^2>(2b)^2-2b^2=2b^2$. Consequently, $p-2b^2>0$ and $p+(c+d)^2-2b^2>0$. This last inequality can be taken to (3.29) to conclude again that $3p > (a+b)^2 - 2d^2$.

Taking the inequality $3p > (a+b)^2 - 2d^2$ to equation (3.28) translates into

$$b(c+d) - (a+b)d < \left(1/\sqrt{2}\right)\sqrt{2(a+b)^2 d^2 + 3p^2 - p^2} - (a+b)d$$

(3.30) $$b(c+d) - (a+b)d < \sqrt{(a+b)^2 d^2 + p^2} - (a+b)d < p\,.$$

The analysis described by (3.27)–(3.30) can be repeated by deducing, also from (3.17), that

(3.31) $$d = +\sqrt{\left[(c+d)^2 - p\right]/2} \qquad (a+b) = +\sqrt{p+2b^2}$$

and concluding, in a similar fashion, that

(3.32) $$(a+b)d - b(c+d) < \sqrt{b^2(c+d)^2 + p^2} - b(c+d) < p\,.$$

Inequalities (3.30) and (3.32) combine to derive that $\left|b(c+d) - (a+b)d\right| < p$. When such a conclusion is brought together with the starting condition, $b(c+d) - (a+b)d \equiv 0 \pmod{p}$, the implication is that $b(c+d) - (a+b)d = 0$, i.e. $b(c+d) = (a+b)d$. Since $\gcd(a+b, b) = 1$ and $\gcd(c+d, d) = 1$, they force $(a+b)\,|\,(c+d)$; say $(c+d) = k(a+b)$, $k \in \mathbb{N}$. Hence $b(c+d) = bk(a+b) = (a+b)d$, which reduces to $d = kb$. But

(3.33) $$p = (a+b)^2 - 2b^2 = (c+d)^2 - 2d^2 = k^2\left[(a+b)^2 - 2b^2\right]$$

implies that $k = 1$. The conclusion is that $b = d$ and $a = c$. Consequently, the two representations (3.17) are identical. □

Lemma 3.2 can also be used to establish that, when two positive integers (prime or not) are represented by binary quadratic forms $F(a,b) = (a+b)^2 - 2b^2$, where $a > b > 0$, their product will be represented likewise, while satisfying the same condition.

**Theorem 3.4.** *Consider two odd integers, p and q (prime or not), that are properly represented by the quadratic form* $F(a,b) = (a+b)^2 - 2b^2$*, where* $a \not\equiv b \pmod 2$ *and* $\gcd(a,b) = 1$*, while satisfying* $a > b > 0$*,*

$$p = (a+b)^2 - 2b^2 \qquad\qquad q = (c+d)^2 - 2d^2\,.$$

*Their product* $N = p \cdot q$ *will also be represented by the form* $F(a,b) = (a+b)^2 - 2b^2$*, with the condition* $a > b > 0$ *equally satisfied.*

*Proof.* From Lemma 3.2, equations (3.13) and (3.14) can be taken to produce

$$N = p \cdot q = (A+B)^2 - 2B^2 = (C+D)^2 - 2D^2$$

$$A = ac + bd, \quad B = b(c+d) + d(a+b), \quad C = a(c+d) + d(a-b), \quad D = bc - ad\,.$$

Since $a > b > 0$ and $c > d > 0$, it can immediately be concluded that $A > 0$, $B > 0$ and $C > 0$, although the sign of $D = bc - ad$ is uncertain.

It can also be deduced that

$$C = a(c+d) + d(a-b) = ac + 2ad - bd > ac - bd > bc - bd > bc - ad = D\,,$$

therefore, $C>D$. The sign of $D$, on the other hand, conveys three possibilities:

(1) $D=0$. In this case, $N=C^2$ is a perfect square. Also $bc=ad$ which, since $\gcd(a,b)=1$ and $\gcd(c,d)=1$, implies that $p=q$.

(2) $D>0$. This satisfies the theorem since $C>D>0$, as required.

(3) $D<0$. Consider the opposite value $D'=-D>0$. It can be written that

$$N=(C+D)^2-2D^2=(C'+D')^2-2D'^2=(C'+D')^2-2D^2\,, \tag{3.34}$$

which implies $(C+D)^2=(C'+D')^2$. It is also the case that

$$D'=ad-bc<ad-bd<ac-bd<ac+2ad-bd=C\,. \tag{3.35}$$

Consequently, $C+D=C-D'>0$, so it is concluded that $C-D'=|C'+D'|>0$. Seeking a positive value for $C'$, it can be assumed that $C'+D'>0$, hence $C-D'=C'+D'$, and finally $C'=C-2D'$

$$C'=a(c+d)+d(a-b)-2\,(ad-bc)=ac-bd+2bc=c(a+b)+b(c-d)>0\,. \tag{3.36}$$

The analysis of the relative values of $C'$ and $D'$ shows that

$$C'=ac+2bc-bd>ad+2bc-bd>ad+2bc-bc>ad-bc>0\,. \tag{3.37}$$

therefore, $C'>D'$, satisfying the theorem.

With regard to the relative values of $A$ and $B$, there are also three possibilities to considered:

(1) $A=B$. This will imply that $N=2A^2$, impossible since $N\equiv 1 \pmod 2$.

(2) $A>B$. This satisfies the theorem since $A>B>0$, as required.

(3) $A<B$. As in the proof of Theorem 3.1, consider a new representation where $A'=A$ and $B'=2A-B$. In such case,

$$\begin{gathered} A'=A=ac+bd>0 \\ B'=2A-B=2ac+2bc-(bc+ad+2bd)=a(c-d)+c(a-b)>0\,. \end{gathered} \tag{3.38}$$

Since $A<B$, $B'=2A-B<A=A'$, concluding that $A'>B'>0$, which satisfies the theorem.

Therefore, the final conclusion is that, in all cases, there will be at least one representation of $N$ by the form $F(a,b)=(a+b)^2-2b^2$ that is compliant with the condition $a>b>0$. □

The application of equation (3.12) to the form $G(a',b')=a'^2-2b'^2$ allows proof by infinite descent of the following: when $\gcd(a',b')=1$, every factor of $N=a'^2-2b'^2$ is itself represented by a form $G(a',b')$ (see Exercise 2.2 in [11]). This can be combined with Theorems 3.1, 3.3 and 3.4 to establish that the sums of the two sides $(x+y)$ of the set of all primitive Pythagorean triangles represent integers $N$ whose factorizations include only primes of the form $p\equiv 1,7 \pmod 8$. Consequently, it is possible to define the set of integers

$$S_{1,7}:=\Big\{N_{1,7}=\prod_{j=1}^{n}p_j^{k_j}=(a+b)^2-2b^2,\ a,b\in\mathbb{N},\ a>b>0,\ \gcd(a,b)=1,\ a\not\equiv b \pmod 2,\ p_j\equiv 1,7 \pmod 8\Big\}\,, \tag{3.39}$$

where $\prod_{j=1}^{n} p_j^{k_j}$ is the prime canonical factorization of $N_{1,7}$.
Since $p_j \equiv 1,7 \pmod 8$, it will always be the case that $N_{1,7} \equiv 1,7 \pmod 8$, explaining the nomenclature used in (3.39).

For all purposes, it can be accepted that the form $F(a,b) = (a+b)^2 - 2b^2$, provided that the condition $a > b > 0$ is satisfied, can be treated as any other positive definite form, despite having a positive discriminant ($\Delta = 8$). This is particularly useful when determining the total number of proper (primitive) representations of $N_{1,7}$ by the form $F(a,b)$. The theory of quadratic forms provides the necessary means to execute the calculation (see Section 8 in chapter VI of [10]). A few considerations are important:

(1) There is only one primitive reduced form with discriminant $\Delta = 8$, this being $G(a',b') = a'^2 - 2b'^2$ (see Table 5 in [27]).

(2) The total number of automorphs, i.e. unimodular substitutions which transform a form into itself [10], for the above primitive reduced form is $w = 2$.

(3) The congruence $h^2 \equiv \Delta \pmod{4N_{1,7}}$, since $\Delta = 8$ and $N_{1,7} \equiv 1 \pmod 2$, hence $\gcd(\Delta, N_{1,7}) = 1$, becomes

$$h^2 \equiv 2 \pmod{N_{1,7}}. \tag{3.40}$$

The total number of incongruent solutions, $h$, of (3.40) (the condition $0 \le h < N_{1,7}$ is required) is the product of the solutions of $h^2 \equiv 2 \pmod{p_j^{k_j}}$, for the various powers $p_j^{k_j}$ composing $N_{1,7}$, as outlined in (3.39). Given that $p_j \equiv 1,7 \pmod 8$, and 2 is a quadratic residue of these primes (for example, see Theorem 11.6 of [28]), it can be concluded that $h^2 \equiv 2 \pmod{p_j}$ always has two incongruent solutions. It can be proved by induction (a direct adaptation of the proof described by Section 5 in Chapter VI of [10]), that the congruences $h^2 \equiv 2 \pmod{p_j^{k_j}}$ still have two solutions when $k_j > 1$. Hence, the total number of solutions of (3.40) is $2^n$ where $n$, as in (3.39), is the number of distinct prime factors of $N_{1,7}$.

(4) The total number of proper representations of $N_{1,7}$ is therefore given by $w$ times $2^n$, this is $2 \cdot 2^n$. These representations fall into groups of 4, which are obtained by changing the signs of $a' = a+b$ and $b' = b$. Since the form $F(a,b) = (a+b)^2 - 2b^2$ has only been defined for positive values of $a$ and $b$, the conclusion is that the actual number of distinct proper representations is $2^{n-1}$.

The properties outlined in this section for the set of integers $S_{1,7}$ (3.39) match closely those described in the previous section for the set $S_{1,5}$ (3.1). Such similarities emphasize the strong geometric parallels between the two binary quadratic forms, and the subsets of prime numbers they encompass.

3.3. $\boldsymbol{x + y - 4r = (a-b)^2 + 2b^2}$. Equation (2.6) indicates that the sum of the two sides of a primitive Pythagorean triangle can be subtracted by four times the radius of the inscribed circle, the result being a binary quadratic form. The new form, $F'(a,b)=(a-b)^2+2b^2$, is found to be equivalent to another form $G'(a',b')$ according to

$$x+y-4r=(a-b)^2+2b^2=a'^2+2b'^2\,, \tag{3.41}$$

(3.41) results from the transformation
$M=\begin{pmatrix}1&0\\1&1\end{pmatrix}$, which induces the change $(a-b,\,b)\to(a,b)$.
The transformation matrix $M$ has determinant $\det(M)=+1$; this makes the new matrix unimodular and justifies the equivalence $F'\sim G'$ [27,30].

The equivalence $F'\sim G'$ determines the type of integers that can be represented by $F'(a,b)=(a-b)^2+2b^2$. As already illustrated by (1.2), all primes $p\equiv 1,3 \pmod 8$ can be expressed by the quadratic form $G'(a',b')=a'^2+2b'^2$ [4,9,11]. It is obvious that, to represent an odd prime, $G'(a',b')$ not only needs $\gcd(a',b')=1$, but also requires $a'\equiv 1 \pmod 2$. The parametrization (2.1) for primitive Pythagorean triangles, where $\gcd(a,b)=1$, $a\not\equiv b \pmod 2$, and $a>b>0$, implies that $a-b>0$, and $a-b\equiv 1 \pmod 2$, as well as $\gcd(a-b,b)=1$. Hence, it can be concluded that the form $F'(a,b)$ can represent odd primes. Furthermore, it is clear that every suitable pair $(a',b')$ can be accommodated by primitive Pythagorean triangles, in line with the equivalence $F'\sim G'$. The application of equation (3.12) to these quadratic forms implies that the previous conclusion can be extended to any integer $N$ whose factorization includes only prime divisors $p\equiv 1,3 \pmod 8$. Therefore, it is possible to define the set of integers

$$\begin{aligned}S_{1,3}:=\Big\{N_{1,3}=\prod p_j^{k_j}=(a-b)^2+2b^2,\, a,b\in\mathbb{N},\, a>b>0,\\ \gcd(a,b)=1,\, a\not\equiv b \pmod 2,\, p_j\equiv 1,3 \pmod 8\Big\}\,,\end{aligned} \tag{3.42}$$

where $\prod p_j^{k_j}$ is the prime canonical factorization of $N_{1,3}$.

Since $p_j\equiv 1,3 \pmod 8$, it is always the case that $N_{1,3}\equiv 1,3 \pmod 8$, and this justifies the nomenclature used in (3.42).

The binary quadratic forms $F'\sim G'$ have a negative discriminant ($\Delta=-8$), thus they are positive definite. To determine the total number of proper (primitive) representations of $N_{1,3}$ by $F'(a,b)=(a-b)^2+2b^2$, the calculation follows steps analogous to those described in the previous section (as before, see Section 8 of Chapter VI in [10]):

(1) There is only one primitive reduced form with discriminant $\Delta=-8$, this being $G'(a',b')=a'^2+2b'^2$ (see table 4 in [27]).

(2) The total number of automorphs for the above primitive reduced form is $w=2$.

(3) The congruence $h^2\equiv\Delta \pmod{4N_{1,3}}$, since $\Delta=-8$ and $N_{1,3}\equiv 1 \pmod 2$ and hence $\gcd(\Delta,N_{1,3})=1$, becomes

$$h^2\equiv 2 \pmod{N_{1,3}}\,. \tag{3.43}$$

The total number of incongruent solutions, $h$, of (3.43) is again the product of the solutions of $h^2\equiv 2 \pmod{p_j^{k_j}}$, for the various powers $p_j^{k_j}$ composing $N_{1,3}$, as in (3.42). Since $p_j\equiv 1,3 \pmod 8$, and $-2$ is a quadratic residue of these primes (see, for example, Exercise 4 in page 201 of [7]), it can be concluded that $h^2 \equiv 2 (\mathrm{mod}\, p_j)$ always has two incongruent solutions. It can again be proved by induction (as in Section 5 of Chapter VI in [10]), that the congruences $h^2 \equiv 2 (\mathrm{mod}\, p_j^{k_j})$ still have only two solutions when $k_j > 1$. Consequently, the total number of solutions of (3.43) is $2^n$ where $n$, as in (3.42), is the number of distinct prime factors of $N_{1,3}$.

(4) Accordingly, the total number of proper representations of $N_{1,3}$ is given by $w$ times $2^n$, this is $2{\cdot}2^n$. These representations fall into groups of 4, obtained by changing the signs of $a' = a - b$ and $b' = b$. Since the form $F(a,b) = (a-b)^2 + 2b^2$ is only defined for positive values of $a$ and $b$, it is concluded that the actual number of distinct proper representations is $2^{n-1}$. As in the preceding sections, when $n = 1$, hence $N_{1,3}$ simplifies to a single prime $p_j \equiv 1{,}3 (\mathrm{mod}\, 8)$, the representation is unique.

Therefore, the properties outlined for the set of integers $S_{1,3}$ (3.42) closely resemble those previously described for the sets $S_{1,5}$ (3.1) and $S_{1,7}$ (3.39). As the overall result, it can be stated that all odd primes are similarly represented by three binary quadratic forms, and that these representations all have a direct geometric interpretation within primitive Pythagorean triangles (Figure 2).

## 4. The triangular theorem of the primes

The results and proofs detailed in Sections 2 and 3 can be summarised in a single theorem.

**Theorem 4.1 (The triangular theorem of the primes).** *Consider the totality of primitive Pythagorean triples, that can be generated by the parametrization*

$$x = 2ab \qquad y = a^2 - b^2 \qquad z = a^2 + b^2$$

*for integers $a, b \in \mathbb{N}$, such that $a > b > 0$, $\gcd(a,b) = 1$, and $a \not\equiv b (\mathrm{mod}\, 2)$. The resulting Pythagoras triangles incorporate geometrically three binary quadratic forms:*

(1) *The sum of the two sides $x$ and $y$ is an integer $N_{1,7}$ represented by*

$$x + y = N_{1,7} = (a+b)^2 - 2b^2,$$

*such that the factorization of $N_{1,7}$ includes only primes $p \equiv 1{,}7 (\mathrm{mod}\ 8)$, which are uniquely represented by the form*

(2) *The result of subtracting two radii, $r = b \cdot (a-b)$, of the inscribed circle from the sum of the two sides $x$ and $y$ is an integer $N_{1,5}$ represented by*

$$x + y - 2r = N_{1,5} = a^2 + b^2 = z,$$

*such that the factorization of $N_{1,5}$ includes only primes $p \equiv 1{,}5 (\mathrm{mod}\ 8)$, which are uniquely represented by the form. $N_{1,5}$ is also the hypotenuse of the triangle, $z$*

(3) *The subtraction of four radii from the sum of the two sides $x$ and $y$ is an integer $N_{1,3}$ represented by*

$$x + y - 4r = N_{1,3} = (a-b)^2 + 2b^2 = z - 2r,$$

*such that the factorization of $N_{1,3}$ includes only primes $p \equiv 1,3 \pmod 8$, which are uniquely represented by the form. The subtraction of two radii from the hypotenuse also yields $N_{1,3}$.*

*Proof.* As detailed in Sections 2 and 3. □

As previously indicated, Figure 2 illustrates the geometrical implications of Theorem 4.1. Key facts are that all odd primes can be geometrically incorporated into primitive Pythagorean triangles and that the three binary quadratic forms have a similar algebraic behaviour, provided the parametrization conditions for $a$ and $b$ are implemented. Table 1 lists some triples and the associated values of $r$, $N_{1,3}$, $N_{1,5}$ and $N_{1,7}$, arising from small values of $a$ and $b$. The primes that occur are highlighted in bold.

TABLE 1. Some primitive Pythagorean triples and associated binary forms.

| $a$ | $b$ | $x$<br>$2ab$ | $y$<br>$a^2-b^2$ | $z$<br>$a^2+b^2$ | $r$<br>$b(a-b)$ | $N_{1,3}$<br>$x+y-4r$ | $N_{1,5}(=z)$<br>$x+y-2r$ | $N_{1,7}$<br>$x+y$ |
|---|---|---|---|---|---|---|---|---|
| 2 | 1 | 4 | 3 | 5 | 1 | **3** | **5** | **7** |
| 3 | 2 | 12 | 5 | 13 | 2 | 9 (=3·3) | **13** | **17** |
| 4 | 1 | 8 | 15 | 17 | 3 | **11** | **17** | **23** |
| 4 | 3 | 24 | 7 | 25 | 3 | **19** | 25 (=5·5) | **31** |
| 5 | 2 | 20 | 21 | 29 | 6 | **17** | **29** | **41** |
| 5 | 4 | 40 | 9 | 41 | 4 | 33 (=3·11) | **41** | 49 (=7·7) |
| 6 | 1 | 12 | 35 | 37 | 5 | 27 (=3·3·3) | **37** | **47** |
| 6 | 5 | 60 | 11 | 61 | 5 | 51 (=3·17) | **61** | **71** |
| 7 | 2 | 28 | 45 | 53 | 10 | 33 (=3·11) | **53** | **73** |
| 7 | 4 | 56 | 33 | 65 | 12 | **41** | 65 (=5·13) | **89** |
| 7 | 6 | 84 | 13 | 85 | 6 | **73** | 85 (=5·17) | **97** |

A further illustration of Theorem 4.1 is provided by Figure 4, where a Venn diagram is shown with three different subsets of odd primes, $P_{1,3}$, $P_{1,5}$ and $P_{1,7}$, classified according to their representation by the three binary quadratic forms outlined in Theorem 4.1:

(1) $P_{1,3} := \{$all primes $p \equiv 1,3 \pmod 8$, represented by $(a-b)^2 + 2b^2$, for some values of $a$ and $b\}$

(2) $P_{1,5} := \{$all primes $p \equiv 1,5 \pmod 8$, represented by $a^2 + b^2$, for some values of $a$ and $b\}$

(3) $P_{1,7} := \{$all primes $p \equiv 1,7 \pmod 8$, represented by $(a+b)^2 - 2b^2$, for some values of $a$ and $b\}$

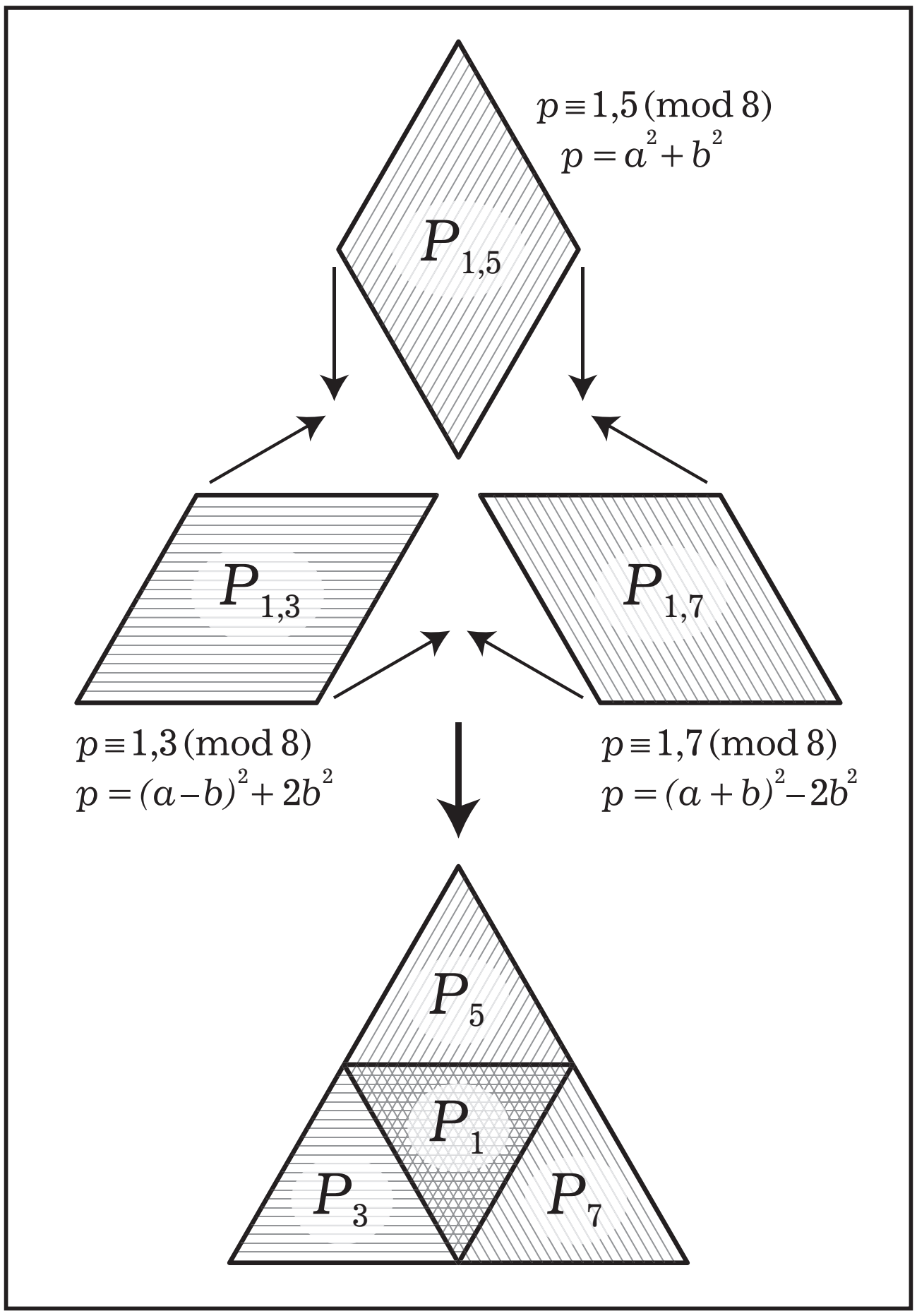

FIGURE 4. The triangular theorem of the primes.

Interestingly, not only the set of all odd primes, $P := \{\text{all primes } p \equiv 1 \pmod 2\}$, results from

$$P = P_{1,3} \cup P_{1,5} \cup P_{1,7}\,, \tag{4.1}$$

but also

$$P_1 = P_{1,3} \cap P_{1,5} \cap P_{1,7} = P_{1,3} \cap P_{1,5} = P_{1,3} \cap P_{1,7} = P_{1,5} \cap P_{1,7}\,, \tag{4.2}$$

where $P_1 := \{\text{all primes } p \equiv 1 \pmod 8\}$ (Figure 4). As also shown in Figure 4, subsets $P_3$, $P_5$ and $P_7$ are defined as containing all primes $p \equiv 3 \pmod 8$, $p \equiv 5 \pmod 8$ and $p \equiv 7 \pmod 8$, respectively.

Since the primes $p \equiv 1 \pmod 8$ can be represented by all three binary quadratic forms, a question to be asked is whether or not a given prime divisor $p \equiv 1 \pmod 8$ can be part simultaneously of the factorization of $N_{1,3}$, $N_{1,5}$ and $N_{1,7}$ for the same primitive triple. The answer is negative: it can be proved that, for any primitive Pythagorean triangle, the three integers $N_{1,3}$, $N_{1,5}$ and $N_{1,7}$ are pairwise relatively prime, i.e. $\gcd(N_{1,3}, N_{1,5}) = 1$, $\gcd(N_{1,3}, N_{1,7}) = 1$ and $\gcd(N_{1,5}, N_{1,7}) = 1$. One proof of this statement makes use of the following lemma.

**Lemma 4.1.** *Consider two odd and distinct integers m and n:* $m \equiv n \equiv 1 \pmod 2$ and $m \neq n$. *If* $\gcd\big((m-n)/2, (m+n)/2\big) = 1$, *then m and n are themselves relatively prime, i.e.* $\gcd(m,n) = 1$.

*Proof.* Contrapositive. First assume that the integers $m$ and $n$ are not relatively prime, i.e. $\gcd(m,n) = d \neq 1$. Hence, it can be written that $m = pd$ and $n = qd$, where $p \equiv q \equiv d \equiv 1 \pmod 2$. Accordingly,

$$m + n = (p+q)\,d \qquad m - n = (p-q)\,d\,, \tag{4.3}$$

where $p \pm q \equiv 0 \pmod 2$. Since $p, q \equiv 1 \pmod 2$, it is elementary to deduce that $p \pm q \equiv 0 \pmod 2$ and $2 \,|\, (p \pm q)$. Consequently,

$$p + q = 2s \qquad p - q = 2t\,, \tag{4.4}$$

and (4.4) can be taken to (4.3) to find

$$(m+n)/2 = sd \qquad (m-n)/2 = td\,. \tag{4.5}$$

Equations (4.5) imply that $d \,\big|\, [(m \pm n)/2]$, i.e. $\gcd\big((m-n)/2, (m+n)/2\big) \neq 1$, hence negating the original condition. Therefore, the final implication is that, if the original condition is satisfied, then $m$ and $n$ are relatively prime. ☐

**Theorem 4.2.** *In any primitive Pythagorean triangle, it will always be found that the integers* $N_{1,7} = (a+b)^2 - 2b^2$, $N_{1,5} = a^2 + b^2$ *and* $N_{1,3} = (a-b)^2 + 2b^2$, *where* $a, b \in \mathbb{N}$, $a > b > 0$, $\gcd(a, b) = 1$, *and* $a \not\equiv b \pmod 2$, *are all pairwise relatively prime, i.e.* $\gcd(N_{1,5}, N_{1,7}) = 1$, $\gcd(N_{1,3}, N_{1,5}) = 1$ *and* $\gcd(N_{1,3}, N_{1,7}) = 1$.

*Proof.* For all three cases, the strategy of the proof is to reach the statement formulated in Lemma 4.1 for the corresponding pairs of integers.

(1) $\gcd(N_{1,5}, N_{1,7}) = 1$.
From the expressions for $N_{1,5}$ and $N_{1,7}$, it can be derived that

$$\begin{aligned} N_{1,7} + N_{1,5} &= (a+b)^2 - 2b^2 + a^2 + b^2 = 2a^2 + 2ab = 2a(a+b) \\ N_{1,7} - N_{1,5} &= (a+b)^2 - 2b^2 - a^2 - b^2 = 2ab - 2b^2 = 2b(a-b)\,, \end{aligned} \tag{4.6}$$

and

$$(N_{1,7} + N_{1,5})/2 = a(a+b) \qquad (N_{1,7} - N_{1,5})/2 = b(a-b)\,. \tag{4.7}$$

Since $\gcd(a,b) = 1$, it is elementary that $\gcd\big(a(a+b), b\big) = \gcd\big(b(a-b), a\big) = 1$. Also, the condition $a \not\equiv b \pmod 2$ combines with $\gcd(a,b) = 1$ to imply that $\gcd(a+b, a-b) = 1$. This can easily be demonstrated with a contrapositive proof: assuming that the opposite is true, i.e. $\gcd(a+b, a-b) \neq 1$, it can then be written that

$$a + b = pd \qquad a - b = qd\,, \tag{4.8}$$

where $d \neq 1$. Given that $a \not\equiv b \pmod 2$, then $a+b \equiv a-b \equiv p \equiv q \equiv d \equiv 1 \pmod 2$ and, accordingly, $p+q \equiv p-q \equiv 0 \pmod 2$. From (4.8), it can be derived that

$$(a+b)+(a-b)=2a=d(p+q)$$
$$(a+b)-(a-b)=2b=d(p-q),$$

and

(4.9) $$a=d(p+q)/2=ds \qquad b=d(p-q)/2=dt.$$

Equations (4.9) imply that $d|a$ and $d|b$, i.e. $\gcd(a,b)\neq 1$, contradicting the starting condition. Hence, $\gcd(a+b,a-b)=1$, as required. From this, it can then be concluded that $\gcd\big(a(a+b),b(a-b)\big)=1$. The last statement can be taken to (4.7) to deduce that $\gcd\big((N_{1,7}+N_{1,5})/2,(N_{1,7}-N_{1,5})/2\big)=1$. Finally, Lemma 4.1 means that $\gcd(N_{1,5},N_{1,7})=1$. □

(2) $\gcd(N_{1,3},N_{1,5})=1$ and $\gcd(N_{1,3},N_{1,7})=1$.
According to (2.2) and (4.6), $N_{1,7}-N_{1,5}=2b(a-b)=2r$, where $r$ is the length of the radius of the inscribed circle for the corresponding primitive Pythagorean triangle. Hence, the preceding conclusion $\gcd(N_{1,5},N_{1,7})=1$ also implies that $\gcd(N_{1,5},r)=\gcd(N_{1,7},r)=1$. Given that $N_{1,5}-N_{1,3}=2r$ and $N_{1,7}-N_{1,3}=4r$, the final implication is that $\gcd(N_{1,3},N_{1,5})=1$ and $\gcd(N_{1,3},N_{1,7})=1$, as required. □

As seen in the proof of Theorem 4.2, the radius $r=b(a-b)$ and the hypothenuse $z=N_{1,5}=a^2+b^2$ of any primitive Pythagorean triangle are relatively prime, i.e. $\gcd(r,z)=1$. The same property is characteristic of the other two binary quadratic forms embedded in the triangle, $N_{1,7}=x+y=(a+b)^2-2b^2$ and $N_{1,3}=x+y-4r=(a-b)^2+2b^2$, i.e. $\gcd(r,N_{1,7})=1$ and $\gcd(r,N_{1,3})=1$. This observation is in sharp contrast with the lack of relative primality most commonly associated with the sides of the triangle, $x$ and $y$: taking from (2.1) that $x=2ab$, it is apparent that $\gcd(r,x)=b\neq 1$, whenever $b>1$; similarly, since $y=a^2-b^2=(a+b)(a-b)$, it is also clear that $\gcd(r,y)=a-b\neq 1$, whenever $a>b+1$. Furthermore, with a single exception (i.e. when $a=2$ and $b=1$, so that $r=1$, $x=4$ and $y=3$), it is evident that no primitive Pythagorean triangle exists where $\gcd(r,x)=1$ and $\gcd(r,y)=1$ occur simultaneously.

## 5. Segregation of primes $p\equiv 1\ (\mathrm{mod}\ 8)$

Using the three binary quadratic forms, as described in Theorem 4.1, it is not possible to separate totally the representation of odd primes according to their congruence modulo 8, given the ubiquitous nature of the primes $p\equiv 1\ (\mathrm{mod}\ 8)$. Nevertheless, the parametrization (2.1) suggests the possibility of differentiating the representations of $N_{1,3}$, $N_{1,5}$ and $N_{1,7}$ according to the specific (and opposing) parities of $a$ and $b$. The following lemma provides the starting result.

**Lemma 5.1.** *The binary quadratic forms* $N_{1,7}=(a+b)^2-2b^2$ *and* $N_{1,3}=(a-b)^2+2b^2$, *where* $a,b\in\mathbb{N}$, $a>b>0$, $\gcd(a,b)=1$, *and* $a\not\equiv b\ (\mathrm{mod}\ 2)$, *always represent integers* $N_{1,7}\equiv N_{1,3}\equiv 1\ (\mathrm{mod}\ 8)$ *when* $a\equiv 1\ (\mathrm{mod}\ 2)$. *Alternatively, when* $a\equiv 0\ (\mathrm{mod}\ 2)$, *the representations are* $N_{1,7}\equiv 7\ (\mathrm{mod}\ 8)$ *and* $N_{1,3}\equiv 3\ (\mathrm{mod}\ 8)$, *respectively.*

*Proof.* The two possibilities for the parities of $a$ and $b$ can be analysed separately.

(1) $a \equiv 1 \pmod 2$, $b \equiv 0 \pmod 2$.
Hence, $a \pm b \equiv 1 \pmod 2$ and $(a \pm b)^2 \equiv 1 \pmod 8$. Also, $b \equiv 0,2,4,6 \pmod 8$, translating into $b^2 \equiv 0,4,0,4 \pmod 8$, respectively, and into $2b^2 \equiv 0 \pmod 8$, for all cases. It can be concluded that $(a \pm b)^2 \mp 2b^2 \equiv 1 + 0 \equiv 1 \pmod 8$. ☐

(2) $a \equiv 0 \pmod 2$, $b \equiv 1 \pmod 2$.
As before, $a \pm b \equiv 1 \pmod 2$ and $(a \pm b)^2 \equiv 1 \pmod 8$. But now $b^2 \equiv 1 \pmod 8$, which implies $\mp 2b^2 \equiv 6,2 \pmod 8$, respectively. The final conclusion is that $(a+b)^2 - 2b^2 \equiv 1 + 6 \equiv 7 \pmod 8$, and $(a-b)^2 + 2b^2 \equiv 1 + 2 \equiv 3 \pmod 8$. ☐

The binary quadratic form $N_{1,5} = a^2 + b^2$ makes no distinction between $a$ and $b$, since both parameters are interchangeable. Nevertheless, when considering the even term of the form, say $b = 2n$, differences can be found depending on the parity of $n$.

**Lemma 5.2.** *The binary quadratic form* $N_{1,5} = a^2 + b^2 = a^2 + 4n^2$, *where* $a, n \in \mathbb{N}$, $\gcd(a, n) = 1$ *and* $a \equiv 1 \pmod 2$, *always represents integers* $N_{1,5} \equiv 5 \pmod 8$ *when* $n \equiv 1 \pmod 2$, *and* $N_{1,5} \equiv 1 \pmod 8$ *when* $n \equiv 0 \pmod 2$.

*Proof.* The two possibilities for the parity of $n$ can be analysed separately.

(1) $n \equiv 1 \pmod 2$.
Since $a \equiv 1 \pmod 2$, it is the case that $a^2 \equiv 1 \pmod 8$. Similarly, $n^2 \equiv 1 \pmod 8$. Therefore, it can be concluded that $a^2 + 4n^2 \equiv 1 + 4 \cdot 1 \equiv 5 \pmod 8$. ☐

(2) $n \equiv 0 \pmod 2$.
Again, $a^2 \equiv 1 \pmod 8$, but now $n^2 \equiv 0 \pmod 8$. Hence, $a^2 + 4n^2 \equiv 1 \pmod 8$. ☐

Using Lemmas 5.1 and 5.2, parity criteria can be taken to the binary quadratic forms representing $N_{1,3}$, $N_{1,5}$ and $N_{1,7}$ (see Theorem 4.1). As a consequence, the resulting expressions will only be able to represent odd primes of a single congruence modulo 8, this is $p \equiv 1,3,5,7 \pmod 8$, respectively.

**Theorem 5.1.** *All odd primes can be represented by binary quadratic forms, according to their congruence modulo* 8*:*

(1) *For all primes* $p \equiv 3 \pmod 8$, $p \equiv 5 \pmod 8$ *and* $p \equiv 7 \pmod 8$, *there will be a single pair of positive integers, s and t, satisfying* $s \equiv t \equiv 1 \pmod 2$ *and* $\gcd(s, t) = 1$, *such that*

(5.1) $p = s^2 + 2t^2$ *when* $p \equiv 3 \pmod 8$,
(5.2) $p = s^2 + 4t^2$ *when* $p \equiv 5 \pmod 8$,
(5.3) $p = s^2 + 4st + 2t^2$ *when* $p \equiv 7 \pmod 8$.

(2) *For all primes* $p \equiv 1 \pmod 8$, *there will be a single pair of positive integers, s and t, satisfying* $s \equiv 1 \pmod 2$, $t \equiv 0,1 \pmod 2$ *and* $\gcd(s, t) = 1$, *such that*

(5.4) $p = s^2 + 8t^2$,
(5.5) $p = s^2 + 16t^2$,
(5.6) *and/or* $p = s^2 + 8st + 8t^2$.

*Proof.* All the representations (5.1)-(5.6) can be derived from Lemmas 5.1 and 5.2.

(1) $p = s^2 + 2t^2$, when $p \equiv 3 \pmod 8$.
Based on Theorem 4.1 and Lemma 5.1, all primes $p \equiv 3 \pmod 8$ can be represented by $p = (a-b)^2 + 2b^2$, where $a > b > 0$, $\gcd(a,b) = 1$, $a \equiv 0 \pmod 2$ and $b \equiv 1 \pmod 2$. Since $a > b$, it is also possible to write $a = s + t$ and $b = t$, so that $s \equiv t \equiv 1 \pmod 2$ and $\gcd(s,t) = 1$. Consequently,

$$p = (a-b)^2 + 2b^2 = (s+t-t)^2 + 2t^2 = s^2 + 2t^2,$$

confirming (5.1).

(2) $p = s^2 + 4t^2$, when $p \equiv 5 \pmod 8$.
Based on Theorem 4.1 and Lemma 5.2, all primes $p \equiv 5 \pmod 8$ can be represented by $p = a^2 + 4n^2$, where $\gcd(a,n) = 1$ and $a \equiv n \equiv 1 \pmod 2$. Thus, making $a = s$ and $n = t$ confirms (5.2).

(3) $p = s^2 + 4st + 2t^2$, when $p \equiv 7 \pmod 8$.
Based on Theorem 4.1 and Lemma 5.1, all primes $p \equiv 7 \pmod 8$ can be represented by $p = (a+b)^2 - 2b^2$, where $a > b > 0$, $\gcd(a,b) = 1$, $a \equiv 0 \pmod 2$ and $b \equiv 1 \pmod 2$. Since $a > b$, it can again be written that $a = s + t$ and $b = t$, so that $s \equiv t \equiv 1 \pmod 2$ and $\gcd(s,t) = 1$. Accordingly,

$$p = (a+b)^2 - 2b^2 = (s+t+t)^2 - 2t^2 = s^2 + 4st + 4t^2 - 2t^2 = s^2 + 4st + 2t^2,$$

confirming (5.3).

(4) $p = s^2 + 8t^2$, when $p \equiv 1 \pmod 8$.
Based on Theorem 4.1 and Lemma 5.1, all primes $p \equiv 1 \pmod 8$ can be represented by $p = (a-b)^2 + 2b^2$, where $a > b > 0$, $\gcd(a,b) = 1$, $a \equiv 1 \pmod 2$ and $b \equiv 0 \pmod 2$. Since $a > b$, write $a = s + b$, where $s \equiv 1 \pmod 2$. Making also $b = 2t$, where $t \equiv 0,1 \pmod 2$ and $\gcd(s,t) = 1$, it is concluded that

$$p = (a-b)^2 + 2b^2 = (s+2t-2t)^2 + 2(2t)^2 = s^2 + 8t^2,$$

confirming (5.4).

(5) $p = s^2 + 16t^2$, when $p \equiv 1 \pmod 8$.
Based on Theorem 4.1 and Lemma 5.2, all primes $p \equiv 1 \pmod 8$ can also be represented by $p = a^2 + 4n^2$, provided that $\gcd(a,n) = 1$, $a \equiv 1 \pmod 2$ and $n \equiv 0 \pmod 2$. Making $a = s$ and $b = 2t$, where $t \equiv 0,1 \pmod 2$ and $\gcd(s,t) = 1$, it is concluded that

$$p = a^2 + 4n^2 = s^2 + 4(2t)^2 = s^2 + 16t^2,$$

confirming (5.5).

(6) $p = s^2 + 8st + 8t^2$, when $p \equiv 1 \pmod 8$.
Based on Theorem 4.1 and Lemma 5.1, all primes $p \equiv 1 \pmod 8$ can be represented by $p = (a+b)^2 - 2b^2$, where $a > b > 0$, $\gcd(a,b) = 1$, $a \equiv 1 \pmod 2$ and $b \equiv 0 \pmod 2$. Since $a > b$, write once more $a = s + b$ where $s \equiv 1 \pmod 2$. Making $b = 2t$, where $t \equiv 0,1 \pmod 2$ and $\gcd(s,t) = 1$, it is concluded that

$$p = (a+b)^2 - 2b^2 = (s+2t+2t)^2 - 2(2t)^2 = s^2 + 8st + 16t^2 - 8t^2 = s^2 + 8st + 8t^2,$$

confirming (5.6).

Equations (5.1)-(5.6) are just restatements of the three binary quadratic forms outlined in Theorem 4.1 that take into account the opposite parities of the positive integers $a$ and $b$, and the condition $a > b$. Therefore, it can be concluded that the representations of all the odd primes by equations (5.1)-(5.6) are unique, exactly as they were for the original quadratic forms. □

To illustrate these results, Table 2 lists the representations of the odd primes $p < 100$ according to the binary quadratic forms (5.1)-(5.6), where all the values of the parameters $s$ and $t$ satisfy the conditions outlined in Theorem 5.1. The table illustrates how these quadratic forms achieve a complete segregation of all odd primes according to their congruence modulo 8, while still providing a unique representation in all cases. The latter fact is especially relevant to the binary forms (5.3) and (5.6), since their discriminants are both positive (8 and 32, respectively) and, consequently, their reduced forms (see Table 5 of [27]), being indefinite, generate an infinite number of representations for primes $p \equiv 3 \pmod 8$ or $p \equiv 1 \pmod 8$, respectively. As already described in Section 3 and Theorem 5.1, the restrictions imposed on the parameters $a$ and $b$ ($s$ and $t$) are the factors that make these representations unique.

Figure 5 further illustrates Theorem 5.1. A Venn diagram can be shown with four different subsets of odd primes, $P_1$, $P_3$, $P_5$ and $P_7$, classified according to their representation by the six binary quadratic forms outlined in Theorem 5.1:

TABLE 2. Segregated representation of primes by binary quadratic forms.

| $p \equiv 3 \pmod 8$<br>$s^2 + 2t^2$ | $p \equiv 5 \pmod 8$<br>$s^2 + 4t^2$ | $p \equiv 7 \pmod 8$<br>$s^2 + 4st + 2t^2$ |
|---|---|---|
| $3 = \mathbf{1}^2 + 2 \cdot \mathbf{1}^2$ | $5 = \mathbf{1}^2 + 4 \cdot \mathbf{1}^2$ | $7 = \mathbf{1}^2 + 4 \cdot \mathbf{1} \cdot \mathbf{1} + 2 \cdot \mathbf{1}^2$ |
| $11 = \mathbf{3}^2 + 2 \cdot \mathbf{1}^2$ | $13 = \mathbf{3}^2 + 4 \cdot \mathbf{1}^2$ | $23 = \mathbf{3}^2 + 4 \cdot \mathbf{3} \cdot \mathbf{1} + 2 \cdot \mathbf{1}^2$ |
| $19 = \mathbf{1}^2 + 2 \cdot \mathbf{3}^2$ | $29 = \mathbf{5}^2 + 4 \cdot \mathbf{1}^2$ | $31 = \mathbf{1}^2 + 4 \cdot \mathbf{1} \cdot \mathbf{3} + 2 \cdot \mathbf{3}^2$ |
| $43 = \mathbf{5}^2 + 2 \cdot \mathbf{3}^2$ | $37 = \mathbf{1}^2 + 4 \cdot \mathbf{3}^2$ | $47 = \mathbf{5}^2 + 4 \cdot \mathbf{5} \cdot \mathbf{1} + 2 \cdot \mathbf{1}^2$ |
| $59 = \mathbf{3}^2 + 2 \cdot \mathbf{5}^2$ | $53 = \mathbf{7}^2 + 4 \cdot \mathbf{1}^2$ | $71 = \mathbf{1}^2 + 4 \cdot \mathbf{1} \cdot \mathbf{5} + 2 \cdot \mathbf{5}^2$ |
| $67 = \mathbf{7}^2 + 2 \cdot \mathbf{3}^2$ | $61 = \mathbf{5}^2 + 4 \cdot \mathbf{3}^2$ | $79 = \mathbf{7}^2 + 4 \cdot \mathbf{7} \cdot \mathbf{1} + 2 \cdot \mathbf{1}^2$ |
| $83 = \mathbf{9}^2 + 2 \cdot \mathbf{1}^2$ | | |
| $p \equiv 1 \pmod 8$<br>$s^2 + 8t^2$ | $p \equiv 1 \pmod 8$<br>$s^2 + 16t^2$ | $p \equiv 1 \pmod 8$<br>$s^2 + 8st + 8t^2$ |
| $17 = \mathbf{3}^2 + 8 \cdot \mathbf{1}^2$ | $17 = \mathbf{1}^2 + 16 \cdot \mathbf{1}^2$ | $17 = \mathbf{1}^2 + 8 \cdot \mathbf{1} \cdot \mathbf{1} + 8 \cdot \mathbf{1}^2$ |
| $41 = \mathbf{3}^2 + 8 \cdot \mathbf{2}^2$ | $41 = \mathbf{5}^2 + 16 \cdot \mathbf{1}^2$ | $41 = \mathbf{3}^2 + 8 \cdot \mathbf{3} \cdot \mathbf{1} + 8 \cdot \mathbf{1}^2$ |
| $73 = \mathbf{1}^2 + 8 \cdot \mathbf{3}^2$ | $73 = \mathbf{3}^2 + 16 \cdot \mathbf{2}^2$ | $73 = \mathbf{5}^2 + 8 \cdot \mathbf{5} \cdot \mathbf{1} + 8 \cdot \mathbf{1}^2$ |
| $89 = \mathbf{9}^2 + 8 \cdot \mathbf{1}^2$ | $89 = \mathbf{5}^2 + 16 \cdot \mathbf{2}^2$ | $89 = \mathbf{3}^2 + 8 \cdot \mathbf{3} \cdot \mathbf{2} + 8 \cdot \mathbf{2}^2$ |
| $97 = \mathbf{5}^2 + 8 \cdot \mathbf{3}^2$ | $97 = \mathbf{9}^2 + 16 \cdot \mathbf{1}^2$ | $97 = \mathbf{1}^2 + 8 \cdot \mathbf{1} \cdot \mathbf{3} + 8 \cdot \mathbf{3}^2$ |

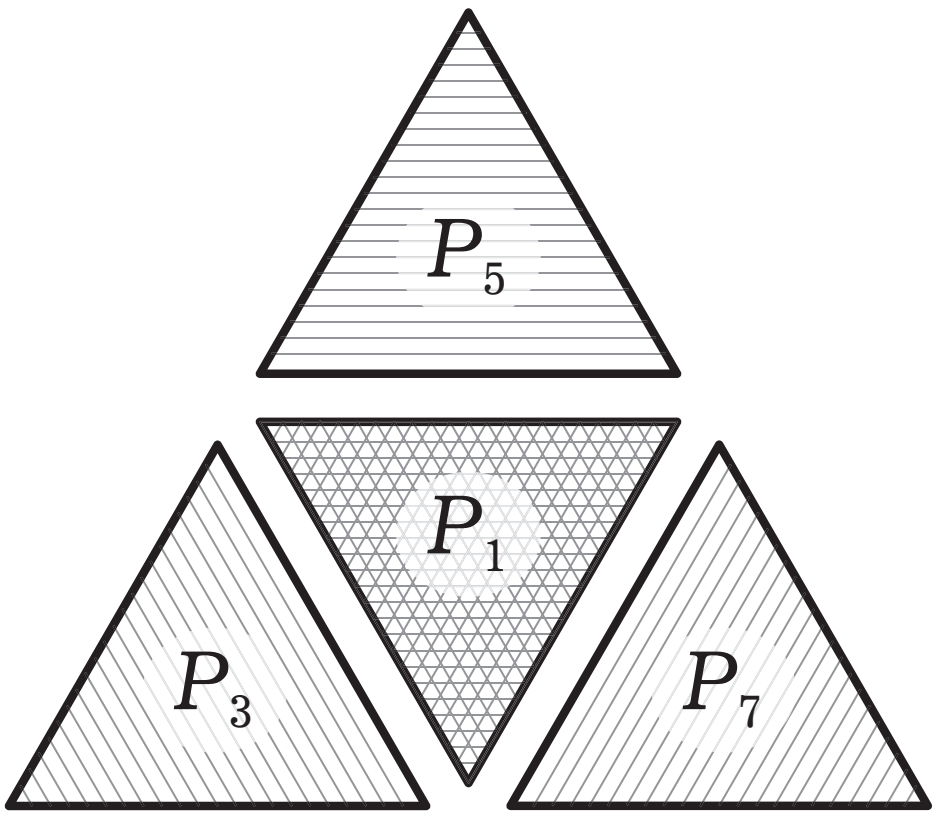

$p \equiv 1 \pmod 8$ $\quad p = s^2 + 8t^2,\ p = s^2 + 16t^2,$
$\quad p = s^2 + 8st + 8t^2$

$p \equiv 3 \pmod 8$ $\quad p = s^2 + 2t^2$
$p \equiv 5 \pmod 8$ $\quad p = s^2 + 4t^2$
$p \equiv 7 \pmod 8$ $\quad p = s^2 + 4st + 2t^2$

FIGURE 5. Schematic segregation of primes and binary quadratic forms.

(1) $P_1 := \{$all primes $p \equiv 1 \pmod 8$, represented by $s^2 + 8t^2$, $s^2 + 16t^2$ and/or $s^2 + 8st + 8t^2$, for some values of $s$ and $t\}$

(2) $P_3 := \{$all primes $p \equiv 3 \pmod 8$, represented by $s^2 + 2t^2$, for some values of $s$ and $t\}$

(3) $P_5 := \{$all primes $p \equiv 5 \pmod 8$, represented by $s^2 + 4t^2$, for some values of $s$ and $t\}$

(4) $P_7 := \{$all primes $p \equiv 7 \pmod 8$, represented by $s^2 + 4st + 2t^2$, for some values of $s$ and $t\}$

Accordingly, the set of all odd primes, $P := \{\text{all primes } p \equiv 1 \pmod 2\}$, results from

(5.7) $$P = P_1 \cup P_3 \cup P_5 \cup P_7 \ .$$

Although the prime segregation described by Theorem 5.1 might be considered superfluous, it is relevant to point out that compliance with the binary quadratic forms included in the theorem extends beyond the prime numbers represented by these forms to all positive integers with factorization including only primes belonging to a single set $P_1$, $P_3$, $P_5$ or $P_7$. This is simple to verify since formula (3.12) can be applied to all six quadratic forms (5.1)-(5.6): the binary forms (5.1), (5.2), (5.4) and (5.5) can be inserted directly into formula (3.12) by making $k = 2, 4, 8$ and 16, respectively; for the quadratic form (5.3), it suffices to use its reduced form $s^2+4st+2t^2=a'^2-2b'^2$ (see Theorem 5.1, as well as Table 5 of [27]), which can be taken to formula (3.12) by making $k=-2$; in the case of the form (5.6), the reduced form is $s^2+8st+8t^2=a'^2-8b'^2$ (see again Theorem 5.1, as well as Table 5 of [27]), and the application of formula (3.12) requires $k=-8$.

The elementary fact that $3\cdot 3\equiv 1 \pmod 8$, $5\cdot 5\equiv 1 \pmod 8$ and $7\cdot 7\equiv 1 \pmod 8$, whereas $3\cdot 1\equiv 3 \pmod 8$, $5\cdot 1\equiv 5 \pmod 8$ and $7\cdot 1\equiv 7 \pmod 8$, respectively, implies that the parity conditions specified in Theorem 5.1 need to be modified to accommodate all possible integers with only one type of prime number(s) in their prime canonical factorization. In a similar fashion to the sets of integers $S_{1,5}$, $S_{1,7}$ and $S_{1,3}$, defined by (3.1), (3.39) and (3.42), respectively, it is also possible to define four distinct sets $S_1$, $S_3$, $S_5$ and $S_7$ as follows:

(1) $S_1 := \Big\{ N_1 = \prod p_j^{k_j} = s^2+8t^2,\ s^2+16t^2$ and/or $s^2+8st+8t^2,\ s,t\in\mathbb{N},\ s,t>0,$
$\gcd(s,t)=1,\ s\equiv 1 \pmod 2,\ t\equiv 0,1 \pmod 2,\ p_j\equiv 1 \pmod 8 \Big\}$,

where $\prod p_j^{k_j}$ is the prime canonical factorization of $N_1$; it is always the case that $N_1\equiv 1 \pmod 8$.

(2) $$S_3 := \Big\{ N_3 = \prod p_j^{k_j} = s^2+2t^2,\ s,t\in\mathbb{N},\ s,t>0,\ \gcd(s,t)=1,\ s\equiv 1 \pmod 2,\ t\equiv \sum k_j \pmod 2,\ p_j\equiv 3 \pmod 8 \Big\},$$

where $\prod p_j^{k_j}$ is the prime canonical factorization of $N_3$; the integers $N_3$ are $N_3\equiv 3 \pmod 8$ when $\sum k_j\equiv 1 \pmod 2$, but $N_3\equiv 1 \pmod 8$ if $\sum k_j\equiv 0 \pmod 2$.

(3) $$S_5 := \Big\{ N_5 = \prod p_j^{k_j} = s^2+4t^2,\ s,t\in\mathbb{N},\ s,t>0,\ \gcd(s,t)=1,\ s\equiv 1 \pmod 2,\ t\equiv \sum k_j \pmod 2,\ p_j\equiv 5 \pmod 8 \Big\},$$

where $\prod p_j^{k_j}$ is the prime canonical factorization of $N_5$; the integers $N_5$ are $N_5\equiv 5 \pmod 8$ when $\sum k_j\equiv 1 \pmod 2$, but $N_5\equiv 1 \pmod 8$ if $\sum k_j\equiv 0 \pmod 2$.

(4) $$S_7 := \Big\{ N_7 = \prod p_j^{k_j} = s^2+8st+8t^2,\ s,t\in\mathbb{N},\ s,t>0,\ \gcd(s,t)=1,\ s\equiv 1 \pmod 2,\ t\equiv \sum k_j \pmod 2,\ p_j\equiv 7 \pmod 8 \Big\},$$

where $\prod p_j^{k_j}$ is the prime canonical factorization of $N_7$; the integers $N_7$ are $N_7\equiv 7 \pmod 8$ when $\sum k_j\equiv 1 \pmod 2$, but $N_7\equiv 1 \pmod 8$ if $\sum k_j\equiv 0 \pmod 2$.

Based on these definitions, it is obvious that the sets $S_1$, $S_3$, $S_5$ and $S_7$ are all subsets of the corresponding sets $S_{1,3}$, $S_{1,5}$ and $S_{1,7}$: while the sets $S_1$ and $S_3$ are both subsets of $S_{1,3}$ (3.42), the sets $S_1$ and $S_5$ are both subsets of $S_{1,5}$ (3.1), and the sets $S_1$ and $S_7$ are both subsets of $S_{1,7}$ (3.39). When determining the total number of proper (primitive) representations for any integer $N_1$, $N_3$, $N_5$ and $N_7$, the latter observation combines with the fact that all such representations are uniquely related to distinct representations of the corresponding integer $N_{1,3}$, $N_{1,5}$ and/or $N_{1,7}$ according to the relevant binary form used to define $S_{1,3}$, $S_{1,5}$ and/or $S_{1,7}$ (a direct implication of Theorem 5.1). Consequently, as already established in Section 3, the number of proper representations for every integer $N_1$, $N_3$, $N_5$ or $N_7$ is given by $2^{n-1}$, where $n$ is the total number of distinct odd primes that divide the integer.

The binary quadratic forms used to define sets $S_{1,3}$, $S_{1,5}$ and $S_{1,7}$ (Section 3 and Theorem 4.1) possess three fundamental and particularly significant properties:

*i) uniqueness*, i.e. only a single representation can be constructed for any prime number expressed by the binary form; *ii) completeness*, i.e. every prime number of the given type can be represented by the corresponding binary form, without exception; and *iii) exclusivity*, i.e. the only integers represented by the binary form are those with a canonical factorization that only includes primes that are also represented by the same form. Given the less restrictive parity conditions introduced in the definitions of sets $S_1$, $S_3$, $S_5$ and $S_7$, it should be no surprise that the prime segregation achieved by Theorem 5.1 is gained at the expense of one of the above properties: although uniqueness and completeness are still satisfied, it is perfectly clear that exclusivity is no longer maintained, since the segregation of prime numbers (based on their congruences modulo 8) cannot be universally transferred to all composite integers with prime canonical factorizations that include primes $p_j \equiv 1 \pmod 8$ alongside primes $p_j \equiv 3 \pmod 8$, $p_j \equiv 5 \pmod 8$ or $p_j \equiv 7 \pmod 8$, respectively.

## 6. Pythagorean prime triplets and conjectures

Theorem 4.1 brings together the three binary quadratic forms characterized in Section 3 and states their intrinsic association with primitive Pythagorean triangles (Figure 2). Given the composition of the integers $N_{1,3}$, $N_{1,5}$ and $N_{1,7}$ that constitute sets $S_{1,3}$ (3.42), $S_{1,5}$ (3.1) and $S_{1,7}$ (3.39), it is natural to consider the cases where all three integers are prime. Hence, the following definition.

**Definition 6.1 (Pythagorean prime triplets).** *Consider a primitive Pythagorean triangle, generated by the parametrization*

$$x = 2ab \qquad y = a^2 - b^2 \qquad z = a^2 + b^2$$

*for integers* $a, b \in \mathbb{N}$, *such that* $a > b > 0$, $\gcd(a,b) = 1$, *and* $a \not\equiv b \pmod 2$. *Consider also the three (geometrically) associated integers* $N_{1,3} = (a-b)^2 + 2b^2 \equiv 1,3 \pmod 8$, $N_{1,5} = a^2 + b^2 \equiv 1,5 \pmod 8$ *and* $N_{1,7} = (a+b)^2 - 2b^2 \equiv 1,7 \pmod 8$. *The resulting triplet* $(N_{1,3}, N_{1,5}, N_{1,7})$ *will be designated a 'Pythagorean prime triplet' if all three integers are prime.*

Table 3 lists some of the smallest Pythagorean prime triplets and makes use of the symbols $p_{1,3}$, $p_{1,5}$ and $p_{1,7}$ to identify the prime numbers according to the respective binary quadratic forms that represent them. With the single exception of the first triplet in the list (i.e. 3, 5 and 7), it can be observed that the length of the radius is always $r \equiv 0 \pmod 3$. There is a simple reason for this pattern. As can easily be verified, for the integer $N_{1,3} = (a-b)^2 + 2b^2$ to be $N_{1,3} \not\equiv 0 \pmod 3$ (required by any prime number rather than 3), the parameters $a$ and $b$, as well as satisfying the condition $\gcd(a,b) = 1$, have to take values such that $a \equiv 1 \pmod 3$ while $b \equiv 0,1 \pmod 3$, or $a \equiv 2 \pmod 3$ while $b \equiv 0,2 \pmod 3$; in all four cases, the radius $r = b \cdot (a-b)$ can be found to be $r \equiv 0 \pmod 3$. Therefore, and with the exception of $(3, 5, 7)$, for all Pythagorean prime triplets it is a necessary (but not sufficient) condition that $3 \mid r$. This implies that the gap between primes must always be $p_{1,7} - p_{1,5} = p_{1,5} - p_{1,3} = 2r \equiv 0 \pmod 6$.

TABLE 3. Some Pythagorean prime triplets.

| $a$ | $b$ | $r$ <br> $b(a-b)$ | $p_{1,3}$ <br> $(a-b)^2+2b^2$ | $p_{1,5}$ <br> $a^2+b^2$ | $p_{1,7}$ <br> $(a+b)^2-2b^2$ |
|---|---|---|---|---|---|
| 2 | 1 | 1 | **3** | **5** | **7** |
| 4 | 1 | 3 | **11** | **17** | **23** |
| 5 | 2 | 6 | **17** | **29** | **41** |
| 10 | 9 | 9 | **163** | **181** | **199** |
| 8 | 3 | 15 | **43** | **73** | **103** |
| 10 | 3 | 21 | **67** | **109** | **151** |
| 10 | 7 | 21 | **107** | **149** | **191** |
| 25 | 24 | 24 | **1153** | **1201** | **1249** |
| 17 | 2 | 30 | **233** | **293** | **353** |
| 14 | 11 | 33 | **251** | **317** | **383** |
| 43 | 42 | 42 | **3529** | **3613** | **3697** |
| 23 | 20 | 60 | **809** | **929** | **1049** |
| 16 | 9 | 63 | **211** | **337** | **463** |
| 35 | 2 | 66 | **1097** | **1229** | **1361** |
| 28 | 25 | 75 | **1259** | **1409** | **1559** |
| 19 | 10 | 90 | **281** | **461** | **641** |
| 23 | 18 | 90 | **673** | **853** | **1033** |
| 26 | 5 | 105 | **491** | **701** | **911** |
| 22 | 15 | 105 | **499** | **709** | **919** |
| 26 | 21 | 105 | **907** | **1117** | **1327** |

As Table 3 suggests, Pythagorean prime triplets appear to be abundant. Given that the gap $(2r)$ between the primes in the triplets is not restricted to very small values, it seems reasonable to anticipate that the progressive decline in prime density encountered in the number line is not a factor that affects the existence of the above triplets in the way it influences the existence of twin primes, or other prime pairs separated by equally small and fixed gaps [14,15,25]. On the other hand, the diminishing density of primes is still certain to have an impact on the likelihood of finding three equally spaced prime numbers that are associated with the same Pythagorean triangle. However, on balance, it can still be said that the appearance of these prime sets along the number line is likely to continue. The expectation of having infinitely many Pythagorean prime triplets is encapsulated in the following conjecture.

**Conjecture 6.2.** *Consider the set of all possible Pythagorean prime triplets. Such a set is infinite.*

TABLE 4. Some Pythagorean prime triplets with all $p \equiv 1 \pmod 8$.

| $a$ | $b$ | $r$ $b(a-b)$ | $p_1$ $(a-b)^2+2b^2$ | $p_1$ $a^2+b^2$ | $p_1$ $(a+b)^2-2b^2$ |
|---|---|---|---|---|---|
| 25 | 24 | 24 | **1153** | **1201** | **1249** |
| 23 | 20 | 60 | **809** | **929** | **1049** |
| 31 | 4 | 108 | **761** | **977** | **1193** |
| 23 | 12 | 132 | **409** | **673** | **937** |
| 35 | 8 | 216 | **857** | **1289** | **1721** |

TABLE 5. Some Pythagorean prime triplets with all $p \not\equiv 1 \pmod 8$.

| $a$ | $b$ | $r$ $b(a-b)$ | $p_3$ $(a-b)^2+2b^2$ | $p_5$ $a^2+b^2$ | $p_7$ $(a+b)^2-2b^2$ |
|---|---|---|---|---|---|
| 2 | 1 | 1 | **3** | **5** | **7** |
| 10 | 9 | 9 | **163** | **181** | **199** |
| 10 | 3 | 21 | **67** | **109** | **151** |
| 10 | 7 | 21 | **107** | **149** | **191** |
| 14 | 11 | 33 | **251** | **317** | **383** |
| 26 | 5 | 105 | **491** | **701** | **911** |
| 22 | 15 | 105 | **499** | **709** | **919** |
| 26 | 21 | 105 | **907** | **1117** | **1327** |
| 22 | 13 | 117 | **419** | **653** | **887** |
| 34 | 21 | 273 | **1051** | **1597** | **2143** |

The ubiquitous character of the primes $p \equiv 1 \pmod 8$ extensively described here implies that, in any Pythagorean prime triplet, these primes are equally likely to be represented by any one of the three binary quadratic forms (Theorem 5.1). Therefore, consideration should be given to the occurrence of triplets where: *i)* all three prime numbers are $p_1 \equiv 1 \pmod 8$ (see Table 4 for some examples), or *ii)* none are, hence they are $p_3 \equiv 3 \pmod 8$, $p_5 \equiv 5 \pmod 8$ and $p_7 \equiv 7 \pmod 8$, respectively (see Table 5). Although the relative size of these two subsets of Pythagorean prime triplets is clearly reduced when compared with the total set of triplets, the expectation is to have infinitely many triplets of both types. Two more conjectures can therefore be formulated.

**Conjecture 6.3.** *Consider the set of all possible Pythagorean prime triplets, where every prime number is* $p \equiv 1 \pmod 8$. *Such a set is infinite.*

**Conjecture 6.4.** *Consider the set of all possible Pythagorean prime triplets, where every prime number is* $p \not\equiv 1 \pmod 8$. *Such a set is infinite.*

Lemmas 5.1 and 5.2 can be brought to the analysis of Pythagorean prime triplets where every prime number is either $p \equiv 1 \pmod 8$ or $p \not\equiv 1 \pmod 8$. As is easily deduced, to ensure that all three primes in the triplet are $p_1 \equiv 1 \pmod 8$, it is required that $a \equiv 1 \pmod 2 \equiv 1,3 \pmod 4$ (Lemma 5.1) and $b \equiv 0 \pmod 4$ (Lemma 5.2); consequently, $r \equiv 0 \pmod 4$ and, since it is also required that $r \equiv 0 \pmod 3$, the final condition is that $r \equiv 0 \pmod{12}$, as confirmed by the examples included in Table 4. The implication is that the gap between primes has to be $2r \equiv 0 \pmod{24}$.

Alternatively, to ensure that all three primes are $p_3, p_5, p_7 \not\equiv 1 \pmod 8$, the requirement is that $a \equiv 0 \pmod 4$ (Lemmas 5.1 and 5.2) and, concomitantly, that $b \equiv 1 \pmod 2 \equiv 1,3 \pmod 4$ (Lemma 5.1); the resulting implication is $r \equiv 1 \pmod 4$. This condition combines with $r \equiv 0 \pmod 3$ to conclude, by applying the Chinese Remainder Theorem [7,10,28], that the global requirement is $r \equiv 9 \pmod{12}$, as the examples given in Table 5 confirm, with the exception of $(3,5,7)$. The implication is that the gap between primes must be $p_7 - p_5 = p_5 - p_3 = 2r \equiv 18 \pmod{24}$.

## 7. Some concluding remarks

7.1. **Pythagorean triangles.** The ability to generate Pythagorean triples seems to have preceded even the formulation of the Pythagorean Theorem itself [1,19]. Despite such an ancient heritage, various aspects of the mathematics associated with Pythagorean triangles still attract the interest of current practitioners (for a survey of recent reports, see [5,6,12,16,17,23,24,26]). Intriguingly, the generating function of Pythagorean triples has also been found to play a significant role in the dynamics of a four-mode quantum system [31]. Consequently, any new observation associated with the basic arithmetic and/or geometry of primitive Pythagorean triangles may lead to additional findings and connections.

The key result of the present article (Theorem 4.1) focuses, not on the generation of primitive Pythagorean triples, but on some of the quadratic forms intimately linked to their geometry. By doing so, it reveals a deep (and previously undetected) connection between the triangles/triples and the totality of odd prime numbers (see Figures 2 and 4). While the identification of the hypothenuses (by equating to the sum of two squares) with odd integers $N_{1,5}$ (3.1) whose factorization only includes primes $p \equiv 1 \pmod 4$, goes back to the work of Fermat [11,30], the presence of two additional binary quadratic forms able to account for all odd integers $N_{1,7}$ (3.39) or $N_{1,3}$ (3.42), with factorizations that include only primes $p \equiv 1,7 \pmod 8$ or $p \equiv 1,3 \pmod 8$, respectively, had clearly been overlooked.

What makes the above observation particularly significant is that all three quadratic forms share striking algebraic similarities, as encapsulated by the three properties listed in Section 5, i.e. *uniqueness* (a single representation for any given prime), *completeness* (every prime of the given type is represented, withouth exceptions), and *exclusivity* (the only integers represented are those with factorization that includes only primes of the given type). Compliance with all three properties is not a universal characteristic of quadratic forms and can only be achieved when the Pythagorean conditions are met by the forms, i.e. for each form $F(a,b)$, the parametrization is such that $a, b \in \mathbb{N}$, $a > b > 0$, $a \not\equiv b \pmod 2$ and $\gcd(a,b)=1$ (see Sections 3.1, 3.2 and 3.3). For any given integer represented by one of the three forms, the total number of primitive representations satisfies exactly the same criterion, and this is especially relevant to $F(a,b)=(a+b)^2-2b^2$ (Section 3.2) which, without the above conditions, would behave as an indefinite quadratic form [22,27].

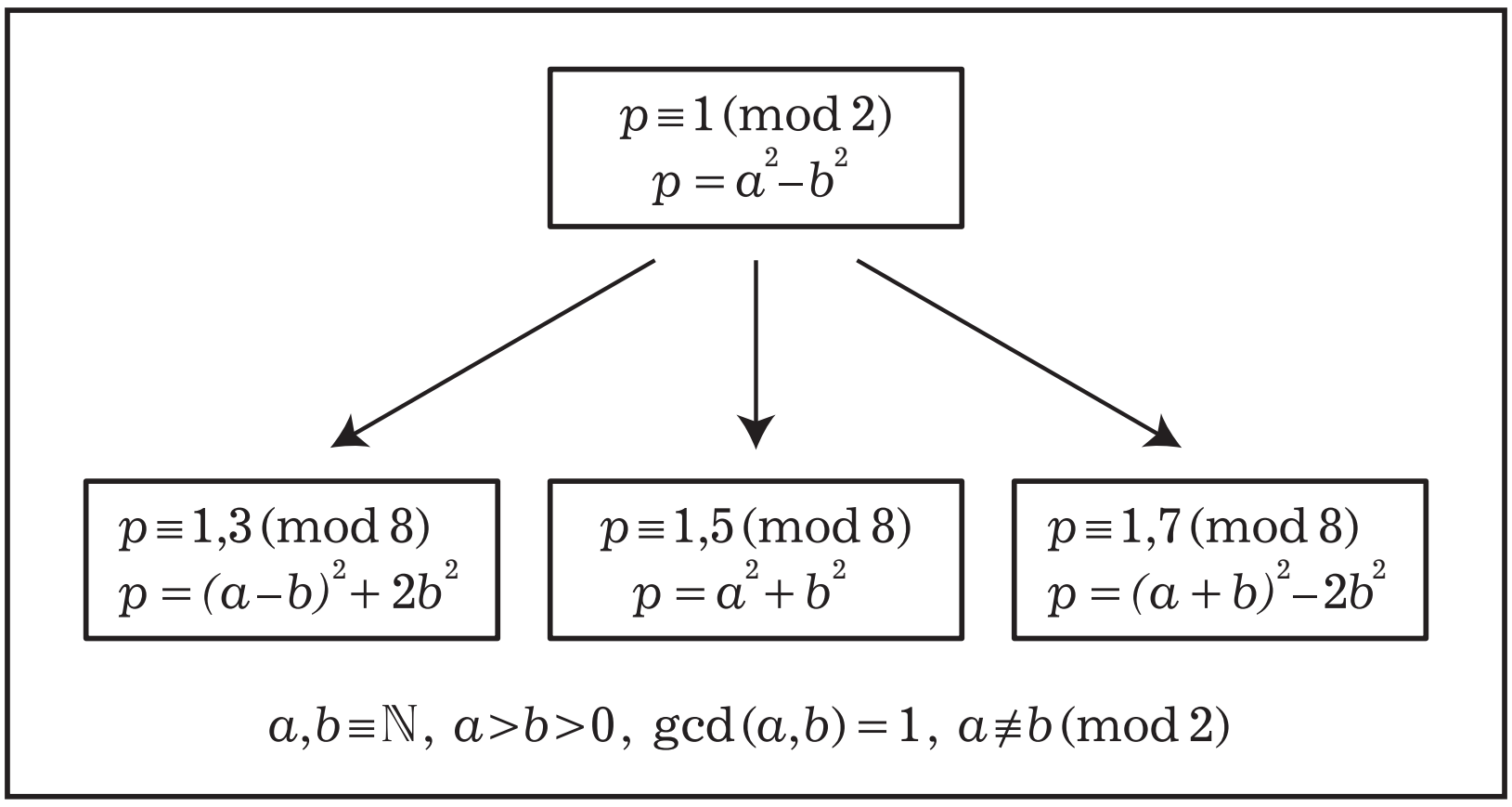

FIGURE 6. Four Pythagorean binary quadratic forms.

Interestingly, these three forms are not the only ones to be found within the geometry of primitive Pythagorean triangles that fully comply with the three properties already discussed. The odd side of the triangles, i.e. $y=a^2-b^2$ (2.1) (Figure 1), accounts for any odd positive integer and the quadratic form involved, i.e. the difference of two squares, also meets all three requirements [22]. Hence, it is elementary that this form accounts for all odd primes $p\equiv 1 \pmod 2$. Figure 6 illustrates how these four Pythagorean binary quadratic forms effectively constitute two layers of representations for all odd prime numbers and their corresponding composite integers, highlighting again the intimate connection between these quadratic forms and the primitive Pythagorean triangles (see Figures 1 and 2). Furthermore, it should be noted that the Pythagorean parametrization (2.1) also serves as a common template for the construction of the four sets of integers and this fact, combined with the algebraic similarities between all four quadratic forms, perhaps might be used to extract other primary results by elementary means.

Another comment can be made on possible ways to classify prime numbers. A distinction is normally made between primes $p\equiv 1 \pmod 4$, represented by the sum of two squares, and primes $p\equiv 3 \pmod 4$, which are not [7,28,29]. This criterion is reinforced by the fact that, when considering Gaussian integers, primes $p\equiv 3 \pmod 4$ are still Gaussian primes, while primes $p\equiv 1 \pmod 4$ are not [28,29]. Although such a classification is reasonable, the two resulting prime sets lack comparable algebraic properties, since there is no single binary quadratic form capable of representing all primes $p\equiv 3 \pmod 4$. The purpose of Figure 4 is to emphasize that Theorem 4.1 provides an alternative classification for all odd primes that makes no formal distinctions between the three resulting sets $P_{1,3}$, $P_{1,5}$ and $P_{1,7}$. Given that subset $P_1$ of primes $p \equiv 1 \pmod 8$ is common to all three prime sets, the result is three interlocking families of prime numbers that (accepting the basically identical algebraic structure of all three sets) produce an aggregate with an apparent three-fold rotational symmetry (Figure 4). Without claiming proof, it seems likely that this combination of binary quadratic forms is unique in representing all odd primes in this manner.

7.2. **Indefinite binary quadratic forms.** Although the theory of quadratic forms is extremely well developed, particularly when concerning positive definite forms [4,9,11,13,21,22,27], their ubiquity in number theory means that further generalizations are possible and the subject of very active research [15,27,29]. Progress is also regularly made in the characterization of indefinite quadratic forms [2,3,8]. Nevertheless, the results reported in this article indicate that findings of a more elementary nature can still be made. In this context, the most significant finding described here is probably that the nature of indefinite binary quadratic forms can be changed under certain simple restrictions (Section 3.2). The case of the form $F(a,b) = (a+b)^2 - 2b^2$ is of immediate interest given the implications for primitive Pythagorean triangles already discussed. But it is also indicative of similar considerations available to a more general case.

When examining primitive reduced forms with positive discriminant $\Delta$ (see, for example Table 5 in [27]), it is clear that, whenever $\Delta \equiv 0 \pmod 4$, there is one indefinite reduced form $G(a',b') = a'^2 - kb'^2$, for certain $k \in \mathbb{N}$, $k > 0$. In line with the transformation used in Section 3.2, it is conceivable that the above form $G(a',b')$ could be transformed into $F(a,b,k,l) = (a+lb)^2 - kb^2$, where $k$ and $l$ are fixed parameters $k, l \in \mathbb{N}$, $k > 0$, $l > 0$, with $a' = a + lb$ and $b' = b$. Hence, the equivalence $F \sim G$. After adopting once more the Pythagorean restrictions on parameters $a$ and $b$, i.e. $a, b \in \mathbb{N}$, $a > b > 0$ and $a \not\equiv b \pmod 2$, in order to ensure that the form $F(a,b,k,l)$ always represents positive integers, i.e. $F(a,b,k,l) = (a+lb)^2 - kb^2 > 0$, it is apparent that a sufficient condition is $(b+lb)^2 \geq kb^2$, translating into $l \geq \sqrt{k} - 1$. Consequently, a suitable selection of the fixed parameters $k$ and $l$ is the only additional prerequisite to avoid negative representations.

Although a systematic evaluation is necessary to establish any potentially useful consequences that might be derived from the above transformations, the results reported in this article for the form $F(a,b) = (a+b)^2 - 2b^2$ suggest that similar findings might be possible, at least in part. As illustration, Table 6 lists some representations for two forms $F(a,b,k,l) = (a+lb)^2 - kb^2$ of particular interest, i.e. $F(a,b,8,3) = (a+3b)^2 - 8b^2$ and $F(a,b,32,5) = (a+5b)^2 - 32b^2$. It is elementary to see that, in both cases, all integers $N$ represented by the form are $N \equiv 1 \pmod 8$. Therefore, any prime numbers represented by these forms will also be $p \equiv 1 \pmod 8$, as Table 6 confirms. Since neither 8 nor 32 are square-free integers, the forms could have been written as $F(a,b,8,3) = (a+3b)^2 - 2(2b)^2$ and $F(a,b,32,5) = (a+5b)^2 - 2(4b)^2$, respectively, both similar to $G(a',b') = a'^2 - 2b'^2$. The implication is that the integers represented by either of the two forms are also represented by $G(a',b') = a'^2 - 2b'^2$. Hence, the only primes that can appear in their factorizations are $p \equiv 1,7 \pmod 8$ (see Section 3.2), as Table 6 also confirms.

TABLE 6. Two binary quadratic forms, no longer indefinite.

| $a$ | $b$ | $(a+3b)^2-8b^2$ | $(a+5b)^2-32b^2$ |
|---|---|---|---|
| 2 | 1 | **17** | **17** |
| 3 | 2 | 49 (=7·7) | **41** |
| 4 | 1 | **41** | 49 (=7·7) |
| 4 | 3 | **97** | **73** |
| 5 | 2 | **89** | **97** |
| 5 | 4 | 161 (=7·23) | **113** |
| 6 | 1 | **73** | **89** |
| 6 | 5 | **241** | 161 (=7·23) |
| 7 | 2 | **137** | 161 (=7·23) |
| 7 | 4 | **233** | 217 (=7·31) |
| 7 | 6 | **337** | 217 (=7·31) |
| 8 | 1 | **113** | **137** |
| 8 | 3 | 217 (=7·31) | **241** |
| 8 | 5 | 329 (=7·43) | 289 (=17·17) |
| 8 | 7 | **449** | **281** |
| 9 | 2 | **193** | **233** |
| 9 | 4 | **313** | 329 (=7·43) |
| 9 | 8 | **577** | **353** |
| 10 | 1 | 161 (=7·23) | **193** |
| 10 | 3 | 289 (=17·17) | **337** |

What makes these two forms interesting is that they appear to represent the same set of primes/integers represented by the form $s^2+8st+8t^2$ (5.6) discussed in Section 5, while still satisfying the Pythagorean parametrization[1]. Therefore, they achieve the same segregation of primes $p \equiv 1 \pmod 8$, without imposing any parity conditions.

A key characteristic, common to all the forms $F(a,b,k,l)=(a+lb)^2-kb^2$ (under the Pythagorean parametric conditions), appears to be that their representations of primes are always unique. Since formula (3.12) can be universally applied to these forms, uniqueness of prime representation will imply that associated composite integers are also representable by the forms, with the total number of primitive representations determined by their canonical factorization as before (Section 3). Although the actual set of prime numbers represented in each case has to be established with care (for example, large values of $k$ and $l$ are certain to preclude the representation of some small primes), uniqueness of representation for primes seems likely enough to prompt the following conjecture.

[1] Not surprisingly, since both $(a+3b)^2-8b^2$ and $s^2+8st+8t^2$ are equivalent to the same reduced form $G(a',b')=a'^2-8b'^2$, and $(a+5b)^2-32b^2=(a+5b)^2-8(2b)^2$.

**Conjecture 7.1.** *Consider the binary quadratic forms* $F(a,b,k,l)=(a+lb)^2-kb^2$ *for integers* $a,b\in\mathbb{N}$, *such that* $a>b>0$, $\gcd(a,b)=1$, *and* $a\not\equiv b\,(\text{mod }2)$, *with fixed parameters* $k,l\in\mathbb{N}$, $k,l>0$ *and* $l\geq\sqrt{k}-1$. *All the representations of prime numbers by these forms are unique.*

As shown here, the representations by the form $F(a,b)=(a+b)^2-2b^2$ exhibit all three properties of uniqueness, completeness and exclusivity. This suggests that the characterization of the more general case might also be of interest, although the property of exclusivity is likely to be particularly elusive.

7.3. ***p* = 2, the lone prime.** It is often claimed that 2 is the oddest of all prime numbers, simply because 2 is the only even prime. While this is an elementary and obvious example of circular reasoning [2], the total absence (in primitive Pythagorean triangles) of representations for 2 comparable to those found for all the other primes can, in fact, be used as a sufficient argument to differentiate this number from all its peers.

Nevertheless, the three binary quadratic forms listed in Theorem 4.1 also provide a criterion for the association of 2 with another type of prime numbers, more specifically with the primes $p\equiv1\,(\text{mod }8)$. By assigning the same value to the two parametrization parameters $a$ and $b$, $a=b$, the output of all three forms becomes the same, i.e. $(a+a)^2-2a^2=a^2+a^2=(a-a)^2+2a^2=2a^2$. And, when $a=1$, all three quadratic forms represent the same number, 2. Therefore, this prime can be represented by all three binary quadratic forms, in common with all primes $p\equiv1\,(\text{mod }8)$, [3] though in sharp contrast with the remaining primes $p\equiv3\,(\text{mod }8)$, $p\equiv5\,(\text{mod }8)$ and $p\equiv7\,(\text{mod }8)$, that are only represented by one of the forms.

While these binary quadratic forms provide a basis to differentiate/associate 2 from/with the other prime numbers, in turn, 2 itself offers another example of the algebraic similarities between the forms. Consider the multiplication of any of the integers $N_{1,3}$, $N_{1,5}$ or $N_{1,7}$ (Theorem 4.1) with a given power of 2, say $2^t$, where $t\in\mathbb{N}$ and $t>0$. As already known for the sums of two squares [22], it is apparent that, in all three cases, primitive representations of the products $2^t\cdot N_{1,3}$, $2^t\cdot N_{1,5}$ or $2^t\cdot N_{1,7}$ (using the corresponding forms) can only be obtained if $t=1$. Applying formula (3.12), and considering one from any of all possible primitive representations of $N_{1,3}$, $N_{1,5}$ or $N_{1,7}$, respectively, the three cases can be analyzed separately:

(1) $2\cdot N_{1,5}=(1^2+1^2)(a^2+b^2)=(1a\mp1b)^2+(1a\pm1b)^2=(a+b)^2+(a-b)^2$.

As shown in the proof of Theorem 4.2, since $a\not\equiv b\,(\text{mod }2)$ and $\gcd(a,b)=1$,

---

[2] The definition of prime number is sufficient to expose this myth. Clearly, 2 is the only prime divisible by 2, in the same way that 3 is the only prime divisible by 3, 5 the only prime divisible by 5, and so on. What generates this distinction between 2 and the other odd primes is nothing more than the use of parities commonly made when characterizing integers.

[3] While this property associates 2 with all primes $p\equiv1\,(\text{mod }8)$, something that still singles out 2 is that the values of the parameters $a$ and $b$ required for its representations by all three forms are the same ($a=b=1$). This is never the case for any of the primes $p\equiv1\,(\text{mod }8)$.

the implication is that $\gcd(a+b, a-b) = 1$, therefore the representation is primitive and unique.

(2) $2 \cdot N_{1,7} = [(1+1)^2 - 2(1)^2][(a+b)^2 - 2b^2] = (A+B)^2 - 2B^2 = (C+D)^2 - 2D^2$.

According to (3.15) (Lemma 3.2), $A = a+b$, $B = a+3b$, $C = a+b$, $D = a-b$. Since $A < B$, as in the proof of Theorem 3.4, $A' = A$ and $B' = 2A - B = a-b$. Consequently, since $\gcd(a+b, a-b) = 1$, $2 \cdot N_{1,7}$ has a single representation that is also primitive.

(3) $2 \cdot N_{1,3} = [(1-1)^2 + 2(1)^2][(a-b)^2 + 2b^2] = (A-B)^2 + 2B^2 = (C-D)^2 + 2D^2$.

The parameters $A$, $B$, $C$ and $D$ can be obtained by applying (3.12) and following steps similar to those used in the proofs of Lemma 3.2 and Theorem 3.4. The result is again that $A = C = a+b$, and $B = D = a-b$. Thus, as in the preceding cases, $2 \cdot N_{1,3}$ has a single primitive representation.

Using the forms $F(a, b) = (a+b)^2 - 2b^2$ and $F'(a, b) = (a-b)^2 + 2b^2$ to construct the primitive representations of $2 \cdot N_{1,7}$ and $2 \cdot N_{1,3}$, respectively, is what makes it possible to determine the similarities between these two cases and $2 \cdot N_{1,5}$. It is simple to verify that the use of formula (3.12) with the corresponding reduced forms, i.e. $G(a', b') = a'^2 - 2b'^2$ and $G'(a', b') = a'^2 + 2b'^2$, would have impeded the detection of these similarities.

JUAN A PEREZ. BERKSHIRE, UK.
*E-mail address*: jap717@juanperezmaths.com